\documentclass[11pt,twoside,a4paper]{article}
\usepackage{graphicx}
\graphicspath{{figsA/}}
\usepackage[font=footnotesize,labelfont=bf]{caption} % to create subfigures
\usepackage[font=footnotesize,labelfont=bf]{subcaption}% to create subfigures

\usepackage[margin=15mm]{geometry}

\usepackage{amssymb}
\usepackage{amsthm}
\usepackage{amsmath}

\numberwithin{equation}{section}

\usepackage{microtype}
\usepackage{mathtools}

\usepackage{bm}
\usepackage{accents}

\newcommand{\RR}{\mathbb{R}}

\newcommand{\SZ}{\mathcal{S}\mathcal{Z}}

\newcommand{\diver}{{\rm div}}%{new command}{old command}
\newcommand{\bnu}{\bm{\nu}}
\newcommand{\bc}{\bm{c}}
\newcommand{\bn}{\bm{n}}
\newcommand{\bx}{\bm{x}}

\newcommand{\bq}{\bm{q}}
\newcommand{\rd}{{\rm d}}

\newtheorem{theorem}{Theorem}
\numberwithin{theorem}{section} % important bit
\newtheorem{lemma}{Lemma}
\numberwithin{lemma}{section} % important bit
\newtheorem{proposition}{Proposition}
\numberwithin{proposition}{section} % important bit
\theoremstyle{definition}\newtheorem{example}{Example}
\numberwithin{example}{section} % important bit
\theoremstyle{definition}\newtheorem{remark}{Remark}
\numberwithin{remark}{section} % important bit
\usepackage{multirow}

\usepackage{listings}
\usepackage{color}

\definecolor{dkgreen}{rgb}{0,0.6,0}
\definecolor{gray}{rgb}{0.5,0.5,0.5}
\definecolor{mauve}{rgb}{0.58,0,0.82}

\providecommand{\keywords}[1]
{
  \small	
  \textbf{\textit{Keywords---}} #1
}

\title{A hybrid finite element method for moving-habitat models in two spatial dimensions}
\makeatletter
\renewcommand\@date{{%
  \vspace{-\baselineskip}%
  \large\centering
  \begin{tabular}{@{}c@{}}
    Jane Shaw MacDonald\textsuperscript{1,} \textsuperscript{2} \\
    \normalsize jane\_macdonald@sfu.ca
  \end{tabular}%
  \quad and\quad
  \begin{tabular}{@{}c@{}}
    Yves Bourgault\textsuperscript{1} \\
    \normalsize ybourg@uottawa.ca
  \end{tabular}%
  \quad and\quad
  \begin{tabular}{@{}c@{}}
    Frithjof Lutscher\textsuperscript{1,} \textsuperscript{3} \\
    \normalsize flutsche@uottawa.ca
  \end{tabular}

  \bigskip

  \textsuperscript{1}Department of Mathematics and Statistics, University of Ottawa, 150 Louis-Pasteur Pvt, Ottawa, K1N 6N5, Ontario, Canada\par
  \textsuperscript{2}Department of Mathematics, Simon Fraser University, 8888 University Dr W, Burnaby, British Columbia, V5A 1S6\par
   \textsuperscript{3}Department of Biology, University of Ottawa, Marie-Curie Private, Ottawa, K1N 9A7, Ontario, Canada\par

  \bigskip

  \date{December 2023}
}}
\makeatother

\begin{document}
\maketitle
\begin{abstract} Moving-habitat models track the density of a population whose suitable habitat shifts as a consequence of climate change. Whereas most previous studies in this area consider 1-dimensional space, we derive and study a spatially 2-dimensional moving-habitat model via reaction-diffusion equations. The population inhabits the whole space. The suitable habitat is a bounded region where population growth is positive; the unbounded complement of its closure is unsuitable with negative growth. The interface between the two habitat types moves, depicting the movement of the suitable habitat poleward. Detailed modelling of individual movement behaviour induces a nonstandard discontinuity in the density across the interface. For the corresponding semi-discretised system we prove well-posedness for a constant shifting velocity before constructing an implicit-explicit hybrid finite element method. In this method, a Lagrange multiplier weakly imposes the jump discontinuity across the interface. For a stationary interface, we derive optimal \textit{a priori} error estimates over a conformal mesh with nonconformal discretisation. We demonstrate with numerical convergence tests that these results hold for the moving interface. Finally, we demonstrate the strength of our hybrid finite element method with two biologically motivated cases, one for a domain with a curved boundary and the other for non-constant shifting velocity. 
\end{abstract}

%\subjclass{65N30, 65M60, 92-10}
%
\maketitle
%%-----------------------------
%%      your text
%%-----------------------------
\section*{Introduction}
\label{sec:intro}
Moving-habitat models are mathematical models that describe the evolution of the density of a biological population when environmental conditions shift in space due to climate change; see for example \cite{Potapov:2004:BullMathBiol,Berestycki:2009:BullMathBiol,MacDonald:2018:JMathBiol,Zhou:2017:JMathBiol,Cobbold:2020:BullMathBiol}. Specifically, as temperature isoclines move towards higher altitude or latitude, the optimal temperature niche of a species moves, so that individuals have to follow in order to stay in optimal conditions. One modelling approach consists of reaction-diffusion equations for the population density over the whole space, divided into two subdomains: a bounded domain where the population growth rate is positive (good habitat) and the unbounded complement of its closure where growth is negative (bad habitat). In contrast to free boundary problems, the spatio-temporal shift of the good habitat is externally prescribed so that the equations are non-autonomous. Matching conditions of flux and density across the interface between good and bad habitat complete the model description, with continuity in population flux but a potential discontinuity in density, representing different movement behaviour and habitat preference of individuals \cite{Maciel:2013:AmNat,Ovaskainen:2003:JApplProbab,MacDonald:2018:JMathBiol}.

Unlike river systems, where an advective term describes the influence of the current on the population density \cite{Lutscher:2011:JTheorBiol,Lutscher:2006:BullMathBiol,Wang:2019:JMathBiol}, a population living in a climate-driven moving habitat does not passively move with the shifting environment. Rather, the advective component of such a system is determined by the shifting properties of the interface between good and bad habitat.

In 1-dimensional space, the good habitat is an interval. Analytical results exist in special cases of constant and periodic shifting speeds when the length of the interval is constant \cite{Potapov:2004:BullMathBiol,Berestycki:2009:BullMathBiol,MacDonald:2018:JMathBiol,Shen:2022:JMathBiol}. We previously developed a finite difference scheme to study the population dynamics for arbitrary movements of the interfaces, for example accelerating shifting speeds and variable length of the good habitat \cite{MacDonald:2021:MathBiosci}.

Realistic landscapes require (at least) 2-dimensional models to capture the heterogeneity that results from topographical variation across space. The simplifying assumption that the landscape is an infinitely wide strip in the direction perpendicular to the shifting direction (as is often done in models for biological invasions) leads to overestimation of the ability of a population to persist in a moving habitat \cite{Phillips:2015:BullMathBiol}. Few analytical studies of moving-habitat models in two dimensions exist \cite{Berestycki:2008:DiscretContinDynS,Berestycki:2009:DiscretContinDynS,Vo:2015:JDifferEquations,Bouhours:2019:JDynDifferEqu}. They typically focus on proving the existence of (forced) traveling waves under simplified geometries, constant shifting speeds, and without consideration of individual movement behaviour. The analytical results in 2-dimensional space are not yet ripe for addressing the questions of concern from ecologists and conservation scientists. To address these concerns, we need a numerical approach that can capture the complex structure of a shifting terrestrial habitat, where even in the simplest case of constant shifting speeds we have still to advance our knowledge on the mechanisms supporting population persistence in 2-dimensional space.

We develop a finite-element method to study a 2-dimensional reaction-diffusion moving-habitat model that considers habitat-dependent dispersal rates and habitat preference. We subdivide space into two distinct domains as indicated above. The interface between these two domains represents the edge of the good habitat. Across the interface, there is a jump in density. It states that the population density approaching the habitat edge from inside the good habitat is proportional to the population density approaching this edge from outside. The porpotionality constant depends on the ratio of the diffusion rates and on the habitat preference \cite{MacDonald:2018:JMathBiol,Maciel:2013:AmNat,Ovaskainen:2003:JApplProbab}. In contrast, the standard jump conditions in domain coupling systems, like those arising in domain decomposition methods and conjugate heat transfer \cite{Quarteroni:1999:domainDecompBook,John:2019:ArchComputMethodE}, describe the difference of the densities across the interface as equal to some measurable function. The challenges in developing the finite-element method for our problem stem from the movement of the interface and the nonstandard jump in density across it.

To resolve the challenge of the moving interface, we employ the arbitrary Lagrangian-Eularian (ALE) description of motion. This method was originally developed for computational fluid dynamics, in particular fluid-structure interaction, where a velocity field is available at every point in the domain at all times \cite{Donea:2006:BookALE}. In our setting, the only points with prescribed velocity lie along the interface. Through an invertible mapping, the ALE description introduces a new coordinate system (`reference frame'), over which the interface is fixed in space. Over this coordinate system, we build our finite element method. With this approach, by exploiting the fact that the location of the interface is known at all times, our interface remains conformal to the mesh. In this way, we avoid complications in implementation due to lack of mesh conformity at the internal interface that can arise in methods like CutFEM and the immersed interface method \cite{Burman:2015:IntJNumerMethEng,Li:2003:TaiwanJMath,Griffith:2020:AnnuRevFluidMech}.

To handle the jump in density across the interface with common finite element methods, the equations describing the jump would appear as an essential condition strongly enforced in the solution space. This approach requires special modifications of the finite element code to handle the lack of conformity of the numerical solution at the interface. Moreover, nonconformal discretisations across the interface can lead to a loss of accuracy of the numerical solution, with a sublinear convergence order. To avoid these difficulties, we propose a new variant of the mortar finite element method; see for example \cite{Bernardi:1994:ANewConform,Belgacem:1999:NumerMath,Bernardi:1993:BookDomDecomp,Mavriplis:1989:phdThesis}. Mortar methods enforce the usual transfer conditions, in particular continuity of the solution at the interface, weakly through an additional variational equation with a dual variable, or Lagrange multiplier. This results in an optimal convergence order, even with nonconformal discretisations, at the expense of solving a saddle point problem.

Our mortar finite element method deviates from the usual ones in two ways, namely in the jump condition itself and in the discretisation of the Lagrange multiplier near vertices of the domain. To our knowledge, only Sacco et al.~enforce a comparable discontinuity of density across a stationary interface, in particular using Lagrange multipliers together with a Three-Field formulation in 3-dimensional space \cite{Sacco:2020:JSciComput}. We introduce the dual variable without the Three-Field formulation to enforce the prescribed discontinuity in density. 

One particular challenge in the analysis of the method is to prove the well-posedness of the variational system, which often rests on the coerciveness of the bilinear form associated to the elliptic operator \cite{Raviart:1977:MathComput,Boffi:2013:BookMixedFEM}. The interface conditions and movement of the interface define a non-coercive bilinear form so that the Lax--Milgram Lemma does not directly apply. The first and most natural mortar method to address our problem leads to a generalised saddle-point problem, for which well-posedness is equivalent to the bilinear forms each satisfying inf-sup conditions \cite{Bernardi:1988:SIAMJNumerAnal,Nicolaides:1982:SIAMJNumerAnal}. In the present paper, we prove well-posedness for important subcases of the general problem with non-zero advective velocity.

We are mainly interested in the accuracy of the primal variable describing the population density. With this consideration, concerning the discretisation of the Lagrange multiplier, we relax the usual condition imposed on its finite element space \cite{Bernardi:1994:ANewConform,Belgacem:1999:NumerMath,Bernardi:1993:BookDomDecomp,Mavriplis:1989:phdThesis}. We define the mortar space from the trace of the finite element space for the primal variable on one side of the interface. This leads to a simple implementation of the method, with no impact on the accuracy and stability for the primal variable. In our analysis, we make no assumption on the mortar space near vertices and still find optimal convergence estimates for the primal variable. This is in contrast with the literature \cite{Bernardi:1994:ANewConform,Belgacem:1999:NumerMath,Bernardi:1993:BookDomDecomp,Mavriplis:1989:phdThesis}.

Finally, our numerical method is an implicit-explicit solver. That is, we apply explicit time-stepping to the nonlinear terms and implicit time-stepping otherwise. The result is an optimally convergent scheme that does not require a nonlinear solver to achieve accurate solutions to a nonlinear system.

This paper is organised as follows: In Section \ref{sec:theModel}, we develop the 2-dimensional analogue of the moving-habitat model presented in \cite{MacDonald:2021:MathBiosci}, that is, we consider arbitrary movement of the interface and a jump in density across this interface. Given an appropriate mapping, we map the physical frame to a reference frame where the suitable habitat is then fixed in space. In Section \ref{sec:semiDiscreteVariationalFormulaitons}, we develop a hybrid variational system for a constant shifting speed of the interface, which is yet easily extendable to more general cases as demonstrated in a later section. The subsequent analysis relies on the assumption on the shifting speed, and we prove that the system is well-posed and equivalent to the original system. In Section \ref{sec:finiteElementFormulation}, assuming a polygonal domain, we develop a finite element method for this system. In Section \ref{sec:errorAnalysis}, in the case of a stationary interface (zero shift), we obtain error estimates for the primal variable. In Section \ref{sec:2dNumerics}, for a moving interface (non-zero shift), we show that the theoretical order of convergence is matched in test cases for linear finite elements. In the same section, we validate our numerical method via comparison to the previously validated 1-dimensional finite difference scheme developed in \cite{MacDonald:2021:MathBiosci}. Finally, we conclude with a demonstration of two biologically motivated applications, that our method can be applied to non-polygonal domains and non-constant shifts of the bounded domain, largely relaxing assumptions made within the analysis.

\section{The Model}
\label{sec:theModel}
In this section, we introduce our system of interest in its most general form. This master equation is posed over a shifting domain, which imposes difficulties for a numerical scheme. We circumvent any need for interpolation due to the shift by a change of variable, which restates the system over a fixed reference frame in which the shifting domain is replaced by an advective term in the system. 

The suitable habitat is represented by an open, bounded region $\Omega_0(t) \subset \mathbb{R}^2,$ dependent on time, $t,$ with a $\mathcal{C}^2$ boundary. We denote its closure by $\bar{\Omega}_0(t)$. The complement of its closure,  $\Omega_1(t) = \mathbb{R}^2\backslash\bar{\Omega}_0(t),$ represents the unsuitable habitat. The interface, $\Gamma(t) = \bar{\Omega}_0(t)\backslash\Omega_0(t),$ represents the habitat edge between the suitable and unsuitable habitat. 

The density of the population over space, $\bx = (x, y),$ and time is denoted by $u(x,y,t)$ and we denote the restriction of $u$ to $\Omega_i$ as 
$$ u|_{\Omega_i} = u_i, \quad i = 0,1.$$
We assume that population dynamics and dispersal happen simultaneously and on the same time scale. Thus, the system representing the population dynamics over space and time is 
\begin{subequations}\label{eq:master2D}
\begin{align}\label{eq:dynamicsOmega02D}
&\partial_tu_0= d_0\Delta u_0 + u_0(r - au_0),  &\text{ in } \Omega_0(t), 
\\ \label{eq:dynamicsOmega12D}
&\partial_tu_1= d_1\Delta u_1 - m_1u_1, &\text{ in } \Omega_1(t), 
\\
\shortintertext{with interface conditions}
\label{eq:DensityCondition2D}
&u_0 = ku_1,   &\text{ on }  \Gamma(t), 
\\ \label{eq:FluxCondition2D}
&d_0\partial_{\bnu}u_0 + (\bc\cdot \bnu)u_0 = d_1\partial_{\bnu}u_1+ (\bc\cdot \bnu)u_1, &\text{ on } \Gamma(t), 
\shortintertext{and asymptotic conditions at infinity}
 \label{eq:limitCondition2D}
&\lim_{|\bx | \to \infty}u(x,y,t) = 0.
\end{align}
\end{subequations}
The operator $\Delta$ is the Laplacian; i.e., $\Delta u = u_{xx} + u_{yy},$ where the subscripts here denote partial derivatives. The coefficients $d_i$ are the constant diffusion rates, $r$ is the intrinsic growth rate, $a$ is the intraspecific competition coefficient, $m_i$ is the mortality rate, $\bc = \bc(x,y,t)$ is the velocity of $\Gamma(t),$ and the proportionality constant across the interface in Equation (\ref{eq:DensityCondition2D}) is 
$$k = k(x,y) = \frac{\alpha(x,y)}{1 - \alpha(x,y)}\sqrt{\frac{d_1}{d_0}}.$$
Here, $\alpha$ represents the probability that an individual enters the suitable habitat when at a point on the interface. The variable $\bnu$ denotes the unit normal pointing outward from $\Omega_0(t)$ into $\Omega_1(t).$ We denote by $\partial_{\bnu}$ the unit normal derivative on $\Gamma(t).$ Equation (\ref{eq:FluxCondition2D}) describes continuity of flux across the interface, we provide its derivation in Appendix \ref{sec:AppA}. System (\ref{eq:master2D}) is the 2-dimensional analogue of the 1-dimensional systems presented in MacDonald and Lutscher \cite{MacDonald:2018:JMathBiol} (moving-habitat model) and Maciel and Lutscher \cite{Maciel:2013:AmNat} (stationary patch) where more details on the formulation, derivation, and biological meaning can be found, particularly with respect to the jump condition.
\subsection{Change of Variable}
To resolve the challenge of the moving interface, we rewrite our system onto a reference frame, wherein we will build our numerical scheme. Detailed calculations are available in \cite{Fernandes:2015:JCompPhys} where the authors derive a mapping for a reference frame of an evolving domain. Here, we provide the derivation of the  system over the reference frame for the bi-domain spatial construction. 

We can group Equations (\ref{eq:dynamicsOmega02D}) - (\ref{eq:dynamicsOmega12D}) and write
\begin{equation}\label{eq:compactlyWritten}
u_{t} = D\Delta u + G(u), \quad \text{ in } \mathbb{R}^2\backslash\Gamma(t),
\end{equation}
where
$$
D = d_i, \quad\text{ on } \Omega_i, \quad\text{and}\quad G(u) = \begin{cases}
u(r - au), \quad\text{ on } \Omega_0,
\\
-mu, \quad\text{ on } \Omega_1.
\end{cases}
$$
We assume, for each $t \in \left[0, T\right],$ that the moving domain, $\bar{\Omega}_0(t),$ is the image of a reference stationary domain  $\bar{\Omega}'_0.$ That is, there exist mappings $x$ and $y$ such that 
$$\bar{\Omega}_0(t) = \left\{(x,y) \in \mathbb{R}^2 \mid x = x(\xi, \eta, t), \ y = y(\xi, \eta, t), \ (\xi, \eta) \in \bar{\Omega}'_0 \right\},$$
where $\bar{\Omega}'_0 \subset \mathbb{R}^2$ is independent of time and $\Omega'_0$ is open and bounded.  We define $\Omega'_1 := \mathbb{R}^2\backslash\bar{\Omega}'_0$ and $\Gamma' := \bar{\Omega}'_0\backslash\Omega_0.$ We assume that all functions involved are smooth enough so that the inverse mappings, $\xi = \xi(x, y, t)$ and $\eta = \eta(x,y,t),$ exist. 

We set $w(\xi(x,y,t), \eta(x,y,t), t) = u(x, y, t),$ where $w = w_i(\xi, \eta, t),$ for $(\xi, \eta)\in \Omega'_i.$ Then, using the chain rule, we compute the first-order temporal partial derivative
\begin{equation}
u_{t} = w_{t} + w_{\xi}\xi_t + w_{\eta}\eta_t = w_{t} + (\mathbf{a}\cdot\nabla_{\xi,\eta})w,
\end{equation}
with $\mathbf{a} = (\xi_t, \eta_t).$ We assume that $w_{\xi\eta} =  w_{\eta\xi},$ so that the second-order spatial partial derivatives are
\begin{align}
u_{xx} &= w_{\xi\xi}(\xi_{x})^2  +  w_{\xi}\xi_{xx} + w_{\eta\eta}(\eta_{x})^2 + w_{\eta}\eta_{xx} +  2w_{\xi\eta}\xi_{x}\eta_{x},
\\
u_{yy} &= w_{\xi\xi}(\xi_{y})^2  +  w_{\xi}\xi_{yy} + w_{\eta\eta}(\eta_{y})^2 + w_{\eta}\eta_{yy} +  2w_{\xi\eta}\xi_{y}\eta_{y}.
\end{align}
Therefore, the dynamics in the $(\xi, \eta, t) - $coordinate system are
\begin{multline}\label{eq:wTransformationGeneral}
w_{t} = D\left(x(\xi,\eta,t),y(\xi,\eta,t)\right)[w_{\xi\xi}\left((\xi_{x})^2 + (\xi_{y})^2\right) + w_{\eta\eta}\left((\eta_{x})^2 +(\eta_{y})^2\right) 
\\
+ w_{\xi}\Delta\xi + w_{\eta}\Delta\eta + 2w_{\xi\eta}\left(\xi_{x}\eta_{x} + \xi_{y}\eta_{y}\right)] - (\mathbf{a}\cdot\nabla_{\xi,\eta})w+ G(w),\quad (\xi, \eta) \in \mathbb{R}^2\backslash\Gamma'.
\end{multline}
For the flux condition, we have that $$\partial_{\bnu}u_i = \left(\bnu \cdot \nabla\right) u,$$
where $\nabla = (\partial_{x}, \partial_{y}).$ Therefore, the flux condition in the $(\xi, \eta, t)$-coordinate system is
\begin{equation}\label{eq:fluxRef}
d_0(\bnu \cdot \nabla)w_0 + w_0\bc\cdot\bnu = d_1(\bnu \cdot \nabla)w_1 + w_1\bc\cdot\bnu, \quad (\xi, \eta) \in \Gamma',
\end{equation}
where 
$$(\bnu\cdot \nabla)w_i = w_{i\xi}\left(\nu_1\xi_{x} + \nu_2\xi_{y}\right) + w_{i\eta}\left(\nu_1\eta_{x} + \nu_2\eta_{y}\right).$$

With these mappings, the proportionality constant in the density-matching condition becomes $\kappa(\xi, \eta, t) = k(x(\xi,\eta,t), y(\xi, \eta, t)).$ 

\begin{example}\label{ex:generalLinearTransform}
We consider an open, bounded domain $\Omega_0(t),$ where the outward unit normal is known on $\Gamma(0).$ If $\Omega_0$ moves only by linear translations $x - c_1t$ and $y - c_2t,$ then $\bc = (c_1, c_2).$ At any point $(x', y')$ on $\Gamma(t)$ the outward unit normal is equal to the outward unit normal at $(x' - c_1t, y' -c_2t)$ on $\Gamma(0)$. Furthermore, the change of variable $\xi = x -c_1t$ and $\eta = y - c_2t$ sets $\Omega_0' = \Omega_0(0)$ and $\Gamma' = \Gamma(0).$ The first-order partial derivatives are all constant and therefore, all second-order partial derivatives vanish. 

Setting $w(\xi, \eta, t) = u(x,y,t),$ we find 
\begin{equation}
w_t = D'(\xi, \eta)\Delta_{\xi,\eta} w + c_1w_\xi + c_2w_\eta + G(w), \quad (\xi, \eta) \in \mathbb{R}^2\backslash\Gamma',
\end{equation}
where $D'(\xi,\eta) := D(\xi + c_1t, \eta + c_2t).$ The flux condition is
\begin{equation}
d_0\left(\nu_1 w_{0\xi} + \nu_2 w_{0\eta}\right) + w_0\left(c_1\nu_1 + c_2\nu_2\right) 
= d_1\left(\nu_1 w_{1\xi} + \nu_2 w_{1\eta}\right) + w_1\left(c_1\nu_1 + c_2\nu_2\right), \quad  (\xi, \eta) \in \Gamma',
\end{equation}
and the density condition is 
$$w_{0} = k(\xi + c_1t, \eta + c_2t)w_{1}, \quad (\xi, \eta) \in \Gamma'.$$
%\qed
\end{example}
Going forward, we set $\xi = x$ and $\eta = y,$ and simplify notation by dropping all prime notation. From now on, we assume linear boundary shifts of the same form as in Example \ref{ex:generalLinearTransform}. We focus our attention to finding solutions to the system
\begin{subequations}\label{eq:2DConstantShift}
\begin{align}
\label{eq:inside2DGeneral}
&\partial_tw = d_i\Delta w + (\bc\cdot\nabla)w + G(w), &\text{ in } \mathbb{R}^2\backslash\Gamma,
\\ \label{eq:densityCondition2DGeneral}
&w_0 = \kappa(x,y)w_1,  &\text{ on }\Gamma,
\\ \label{eq:fluxCondition2DGeneral}
&d_0\partial_{\bnu}w_0 + (\bc\cdot\bnu)w_0 = d_1\partial_{\bnu}w_1 + (\bc\cdot\bnu)w_1,  &\text{ on }\Gamma,
\\ \label{eq:outerBoundary2DGeneral}
&w_1 \to 0, \text{ as } |(x,y)| \to \infty.
\end{align}
\end{subequations}

Within the present analysis, we will make assumptions on the properties of the shifting speed, and at times, for analytical purposes or to define auxiliary problems, take the shifting speed to be zero. In implementation, we resolve these limitations by matching numerical results in more general scenarios of non-zero shifts to our theoretical conclusions. 
\section{Semi-discrete Variational Formulations}\label{sec:semiDiscreteVariationalFormulaitons}
We are mainly interested in the spatial discretisation of our model. Particularly, our goal is to derive a scheme that can treat the interface conditions. Within this section, we assume that both $\bc$ and $\kappa$ are constant. We approximate the time derivative so that $\partial_t w \approx \frac{w^{n+1} - w^n}{\tau},$ where $\tau$ is the time step. We use a forward-backward Euler scheme; i.e., we take all the linear terms at the $(n+1)^{th}$ time step and the nonlinear term at the $n^{th}$ time step. Then the system at every time step becomes 
\begin{subequations}\label{eq:2DSemiDiscrete}
\begin{align}
\label{eq:inside2DSemiDiscrete}
&\frac{w^{n+1}}{\tau} - d_i\Delta w^{n+1} - (\bc\cdot\nabla)w^{n+1} = G(w^n) + \frac{w^n}{\tau}, &\text{ on }\mathbb{R}^2\backslash\Gamma,
\\ \label{eq:densityConditionSemiDiscrete}
&w_0^{n+1} = \kappa w_1^{n+1},  &\text{ on }\Gamma,
\\ \label{eq:fluxConditionSemiDiscrete}
&d_0\partial_{\bnu}w_0^{n+1} + (\bc\cdot\bnu)w_0^{n+1} = d_1\partial_{\bnu}w_1^{n+1} + (\bc\cdot\bnu)w_1^{n+1},  &\text{ on }\Gamma,
\\ \label{eq:outerBoundarySemiDiscrete}
&w_1^{n+1} \to 0, \text{ as } |\bx| \to \infty.
\end{align}
\end{subequations}

For our variational formulation, we take $\Omega$ to be a bounded domain. We let $\Omega = \Omega_0 \cup \Omega_1 \cup \Gamma.$ We assume homogeneous Dirichlet boundary conditions on the boundary of $\Omega,$ which we denote by $\partial\Omega.$ That is, we consider the system
\begin{subequations}\label{eq:2DSemiDiscreteTrunc}
\begin{align}
\label{eq:inside2DSemiDiscreteTrunc}
&\frac{w^{n+1}}{\tau} - d_i\Delta w^{n+1} - (\bc\cdot\nabla)w^{n+1} = G(w^n) + \frac{w^n}{\tau}, &\text{ on }\Omega\backslash\Gamma,
\\ \label{eq:densityConditionSemiDiscreteTrunc}
&w_0^{n+1} = \kappa w_1^{n+1},  &\text{ on }\Gamma,
\\ \label{eq:fluxConditionSemiDiscreteTrunc}
&d_0\partial_{\bnu}w_0^{n+1} + (\bc\cdot\bnu)w_0^{n+1} = d_1\partial_{\bnu}w_1^{n+1} + (\bc\cdot\bnu)w_1^{n+1},  &\text{ on }\Gamma,
\\ \label{eq:outerBoundarySemiDiscreteTrunc}
&w_1= 0, &\text{ on }\partial\Omega. % &d_1\frac{\partial w_1}{\partial \bm{n}} + (c\cdot n)w_1 = 0, &\bx\in\partial\Omega.
\end{align}
\end{subequations}

With this formulation, we are no longer studying a nonlinear problem, as the left-hand side containing the nonlinearity of the problem is known at every time step. Going forward, we replace $w^{n+1}$ with $w$ and $G(w_i^n) + \frac{w_i^n}{\tau}$ with $g_i.$ Our proofs hold for any $g_i \in L^2(\Omega_i).$

Before we introduce our hybrid formulation in which we can weakly enforce the jump condition across $\Gamma,$ we construct some auxiliary bilinear variational formulations, which will assist us in proving the well-posedness of our desired hybrid formulation. In these bilinear formulations, the jump across the interface is an imposed condition on the solution space. These variational formulations are introduced in Section \ref{sec:bilinearVariationalFormulation}. Later in Section \ref{sec:hybridVariationalForm} we introduce the hybrid formulation for which our numerical method is built.
\subsection{Bilinear Variational Formulation}\label{sec:bilinearVariationalFormulation}
Let 
$$X = \{v\in L^2(\Omega) : v|_{\Omega_i} \in H_*^1(\Omega_i)\} \equiv \Pi_{i = 0}^1H_*^1(\Omega_i),$$
where $H_*^1(\Omega_0) = H^1(\Omega_0)$ and $H_*^1(\Omega_1) = H_{\partial\Omega, 0}^1(\Omega_1) = \{ v \in H^1(\Omega_1) : v = 0 \text{ on } \partial\Omega\}.$ We endow the space $X$ with the norm 
$$\left\|v\right\|_{X} = \left(\sum_{i = 0}^1\left\|u_i\right\|_{1, \Omega_i}^2\right)^\frac{1}{2},$$
where $\left\| \cdot \right\|_{1, \Omega_i}$ denotes the typical $H^1$ norm over $\Omega_i.$

For a function $w \in X,$ we denote the restriction of $w$ to $\Omega_i$ as $w|_{\Omega_i} = w_i.$
Let 
$$\mathcal{V}_{\kappa} = \{v \in X : v_0 = \kappa v_1 \text{ on } \Gamma \}.$$ 
Then,
$$\mathcal{V}_1 = \{v \in X : v_0 = v_1  \text{ on } \Gamma \} = H_0^1(\Omega).$$
The last equality results from the fact that $v_0 = v_1  \text{ on } \Gamma.$ To obtain a variational formulation, we multiply Equation (\ref{eq:inside2DSemiDiscreteTrunc}) by $v \in \mathcal{V}_1,$ integrate, and find the relation
\begin{equation}\label{eq:noScale}
\sum_{i = 0}^1 \int_{\Omega_i}\left[-d_i\Delta w_i - \bm{c}\cdot\nabla w_i + \frac{w_i}{\tau}\right]v_i \ \rd\bx = \sum_{i = 0}^1 \int_{\Omega_i} g_i v_i  \ \rd\bx.
\end{equation}
We apply Green's formula and integration by parts to the left-hand side of this equation and thus derive our first variational problem, which we call (\textbf{P}).
\begin{description}
\item[($\text{\textbf{P}}$)]\label{prob:bilinearForm1}
Find $w \in \mathcal{V}_{\kappa}$ such that
\begin{equation}\label{eq:variationalP}
a(w,v) = \ell(v), 
\end{equation}
for all $v \in \mathcal{V}_1,$ where
\begin{equation}\label{eq:a(w,v)}
a(w,v) =\sum_{i = 0}^1\int_{\Omega_i}d_i\nabla w_i \nabla v_i + (\bc\cdot\nabla v_i + \frac{v_i}{\tau})w_i\ \rd\bx, 
\end{equation}
and 
\begin{equation}\label{eq:L(v)}
\ell(v) = \sum_{i = 0}^1 \int_{\Omega_i} g_i v_i  \ \rd\bx.
\end{equation}
\end{description}
There are no terms integrating along $\Gamma,$ since $v_1 = v_0$ on $\Gamma$ and $w$ satisfies Equation (\ref{eq:fluxConditionSemiDiscreteTrunc}).

This problem is encumbered with the fact that the test space and the solution space are not the same. In other words, the Lax--Milgram Lemma does not apply. In this case, to directly show that Problem (\textbf{P}) is well-posed, we could try to show that the conditions of the Banach-Ne\u{c}as-Babu\u{s}ka Theorem \cite{Ern:2004:BookTheoryAndPractise} are satisfied. Yet, we circumvent this challenge by finding an equivalent bilinear variational form to which the Lax--Milgram Lemma is applicable. 

Let us consider a new variational problem, which we call ($\text{\textbf{P}}'$).
\begin{description}
\item[($\text{\textbf{P}}'$)]\label{prob:bilinearForm2}
Find $w \in \mathcal{V}_{\kappa}$ such that
\begin{equation}\label{eq:variationalPPrime}
\tilde{a}(w,v) = \tilde{\ell}(v), 
\end{equation}
for all $v \in \mathcal{V}_\kappa,$ where
\begin{equation}\label{eq:aTilde(w,v)}
\begin{split}
\tilde{a}(w,v) = \int_{\Omega_0}d_0\nabla w_0 \nabla v_0 + (\bc\cdot\nabla v_0 + \frac{v_0}{\tau})w_0 \ \rd\bx 
\\ +  \ \kappa\int_{\Omega_1}d_1\nabla w_1 \nabla v_1 + (\bc\cdot\nabla v_1 + \frac{v_1}{\tau})w_1 \ \rd\bx, 
\end{split}
\end{equation}
and 
\begin{equation}\label{eq:LTilde(v)}
\tilde{\ell}(v) =  \int_{\Omega_0} g_0 v_0  \ \rd\bx +  \int_{\Omega_1} \kappa g_1 v_1  \ \rd\bx.
\end{equation}
\end{description}

For Problem ($\text{\textbf{P}}'$), the test functions and the solution belong to the same space. Thus, for well-posedness, we need only show that the conditions of the Lax--Milgram Lemma are satisfied. Moreover, solving Problem ($\text{\textbf{P}}'$) is equivalent to solving Problem (\textbf{P}); see Theorem \ref{thm:PUnique}. 
\subsection{Existence and Uniqueness for Problem ($\text{\textbf{P}}'$).}\label{sec:existAndUniqeProbP'}
We divide the proof of existence and uniqueness for Problem ($\text{\textbf{P}}'$) into two cases, the first, the no shift case, that is $\bc = 0$, the second, the non-zero shift case, $\bc \neq 0$. In the first case, we make no restrictions on the other model parameters. In the second case, we take $d_1 = d_0$ to ensure coercivity of the bilinear form. Such restrictions on the diffusion coefficient are not uncommon in population models, particularly when diffusion rates are difficult to estimate. Additionally, when $d_1 \neq d_0,$ this defines a non coercive bilinear form and proving well-posedness of that problem reaches beyond the scope of this analysis. The proof of existence and uniqueness in these cases of a nonstandard jump across a stationary interface with differing diffusion rates and a nonstandard jump across a moving interface with a globally constant diffusion rate is novel to the theory.

We apply the Lax--Milgram Lemma to prove that Problem ($\text{\textbf{P}}'$) is well-posed. Let $\tilde{a}_0(\cdot, \cdot) \in \mathcal{L}\left(\mathcal{V}_\kappa \times \mathcal{V}_\kappa; \mathbb{R}\right)$ denote the bilinear form $\tilde{a}(\cdot, \cdot)$ when $\bc = 0.$ Let ($\text{\textbf{P}}'_0$) denote the variational problem obtained by setting $\bc = 0$ in Problem ($\text{\textbf{P}}'$); i.e., 
\begin{description}
\item[($\text{\textbf{P}}'_0$)]\label{prob:bilinearForm2_c0}
Find $w \in \mathcal{V}_{\kappa}$ such that
\begin{equation}\label{eq:variationalProblem0}
\tilde{a}_0(w,v) = \tilde{\ell}(v),
\end{equation}
for all $v \in \mathcal{V}_\kappa.$
\end{description}

\begin{theorem}\label{theo:c0}
Problem ($\text{\textbf{P}}'_0$) is well posed. Moreover, when $\bc = 0,$ $w \in \mathcal{V}_{\kappa}$ is the unique solution of System (\ref{eq:2DSemiDiscreteTrunc}) if and only if it is the unique solution of Problem ($\text{\textbf{P}}'_0$) 
\end{theorem}
\begin{proof}
First, note that $\mathcal{V}_{\kappa}$ is a Hilbert space equipped with the norm from $X.$ Indeed, the following mapping is linear and continuous:
$$
\text{
\begin{tabular}{c c c c c c c} 
$X$ &$\to$ &$H^1(\Omega_0)\times H_{0,\partial \Omega}^1(\Omega_1)$ &$\to$ &$H^{1/2}(\Gamma)\times H^{1/2}(\Gamma)$ &$\to$ &$H^{1/2}(\Gamma)$
\\
$v$ &$\mapsto$ &$(v_0, v_1)$ &$\mapsto$ &$(v_0|_{\Gamma}, v_1|_{\Gamma})$ &$\mapsto$ &$v_0 - \kappa v_1.$
\end{tabular} 
}
$$
Let us denote this jump operator by $[v]_{\kappa} =  v_0 - \kappa v_1.$ We have that $\mathcal{V}_{\kappa} = ([v]_{\kappa})^{-1}(\{0\});$  i.e., $\mathcal{V}_{\kappa}$ is the preimage of $\{0\}.$ Since the set $\{0\}$ is closed in every vector space, $\mathcal{V}_{\kappa}$ is a closed subspace of $X.$ 

The bilinear form $\tilde{a}_0(u, v)$ is coercive, $\tilde{a}_0(w, w) \geq  \underline{\sigma} \left\|w\right\|_X,$ and bounded, $|\tilde{a}_0(w, v)| \leq \overline{\sigma} \left\|w\right\|_X^2\left\|v\right\|_X^2,$ where $\underline{\sigma}$ and $\overline{\sigma}$ are dependent only on $d_0, \ d_1, \ \tau,$ and $\kappa.$ 

Thus, from the Lax--Milgram Lemma, there exists a unique $w \in \mathcal{V}_{\kappa}$ such that 
$$\tilde{a}_0(w, v) = \tilde{\ell}(v), \quad \forall v \in \mathcal{V}_{\kappa},$$ 
and
$$\left\|w\right\|_X \leq C\left\|g\right\|_{0,\Omega}, \ \ C > 0 \text{ is a constant.}$$
Here, $\left\| \cdot \right\|_{0, \Omega}$ denotes the $L^2$ norm over $\Omega.$

Now, suppose that $w \in \mathcal{V}_{\kappa}$ is the unique solution to Problem ($\text{\textbf{P}}'_0$). Then, for any $v \in \mathcal{D}(\Omega_0)$ Equation (\ref{eq:variationalProblem0}) is satisfied. Applying Green's formula, we find that
 $$\left\langle-d_0\Delta w_0 + \frac{1}{\tau}w_0, v_0\right\rangle = \int_{\Omega_0} g_0v_0 \ \rd \bx, \qquad \forall v \in \mathcal{D}(\Omega_0).$$
 Similarly, we see that 
  $$\left\langle\kappa\left(-d_1\Delta w_1 + \frac{1}{\tau}w_1\right), v_1 \right\rangle = \int_{\Omega_1}\kappa g_1v_1 \ \rd \bx, \qquad \forall v \in \mathcal{D}(\Omega_1).$$
 Thus, we see that $w$ satisfies Equation (\ref{eq:inside2DSemiDiscreteTrunc}) in the sense of distributions when $\bc = 0.$ We set $\mathfrak{D} = \cup_{i = 0}^1\mathcal{D}(\bar{\Omega}_i).$ Then, taking any $v \in \mathcal{V}_{\kappa}\cap\mathfrak{D},$ and applying Green's formula, we find that 
\begin{multline*}
\left\langle-d_0\Delta w_0 + \frac{1}{\tau}w_0 - g_0, v_0 \right\rangle + \left\langle\kappa\left(-d_1\Delta w_1 + \frac{1}{\tau}w_1 - g_1\right), v_1 \right\rangle 
\\
+ \left\langle d_0\partial_{\bnu}w_0, v_0 \right\rangle_{\Gamma} - \left\langle \kappa d_1\partial_{\bnu}w_1, v_1 \right\rangle_{\Gamma} = 0, \qquad \forall v \in \mathcal{V}_{\kappa}. 
\end{multline*}
Here, $\langle \cdot, \cdot \rangle_{\Gamma}$ is the duality product between $H^{-1/2}(\Gamma)$ and $H^{1/2}(\Gamma).$ Since Equations (\ref{eq:inside2DSemiDiscreteTrunc}) are satsified in the sense of distributions, we conclude that the boundary term must be equal to zero. Since for $v \in \mathcal{V}_{\kappa},$ we have that $v_0 = \kappa v_1,$ this boundary equation becomes 
$$\left\langle d_0\partial_{\bnu}w_0 -  d_1\partial_{\bnu}w_1, \kappa v_1\right\rangle_{\Gamma} = 0, \qquad \forall v \in \mathcal{V}_{\kappa}.$$
Hence, we recover the flux condition (\ref{eq:fluxConditionSemiDiscreteTrunc}), given that $\bc = 0.$

Next, suppose that $w$ satisfies System (\ref{eq:2DSemiDiscreteTrunc}) in the sense of distributions. We multiply the equation inside $\Omega_1$ corresponding to Equation (\ref{eq:inside2DSemiDiscreteTrunc}) by $\kappa.$ We multiply the resulting coupled equations by $v \in \mathcal{V}_{\kappa}\cap\mathfrak{D},$ integrate and then find the relation
\begin{multline}\label{eq:beforeP0'}
\left\langle-d_0\Delta w_0 - (\bc\cdot\nabla) w_0 + \frac{w_0}{\tau}, v_0 \right\rangle + \left\langle\kappa\left(-d_1\Delta w_1 - (\bc\cdot\nabla) w_1 + \frac{w_1}{\tau}\right), v_1\right\rangle
\\
= \int_{\Omega_0}g_0v_0 \ \rd \bx + \int_{\Omega_1}\kappa g_1v_1 \ \rd \bx. 
\end{multline}
Applying Green’s formula and integration by parts transforms Equation (\ref{eq:beforeP0'}) into Problem ($\text{\textbf{P}}'$), considering that $v_0 =\kappa v_1$ on $\Gamma$ for $v\in\mathcal{V}_{\kappa}$ and that $w$ satisfies Equation (\ref{eq:fluxConditionSemiDiscreteTrunc}).
%\qed
\end{proof}
We now consider the case when $\bc \neq 0$ and $d_0 = d_1 = d.$ First, following Chapter 2 of \cite{Cantrell:2004:book} and considering System (\ref{eq:2DSemiDiscreteTrunc}), we make the change of variable $w_i = e^{-\frac{\bc}{d} \cdot \bx}u_i.$ Then, we have that 
$$\nabla w_i = -\frac{\bc}{d} e ^{-\frac{\bc}{d}\cdot \bx}u_i +e^{-\frac{\bc}{d} \cdot \bx}\nabla u_i,$$
and
$$\Delta w_i = \nabla\cdot\nabla w_i = \left|\frac{\bc}{d}\right|^2e^{-\frac{\bc}{d} \cdot \bx}u_i - 2e ^{-\frac{\bc}{d} \cdot \bx}\frac{\bc}{d}\cdot\nabla u_i + e^{-\frac{\bc}{d} \cdot \bx}\Delta u_i.$$
Hence, Equation (\ref{eq:inside2DSemiDiscreteTrunc}) becomes
\begin{equation}\label{eq:notDivForm}
-d  e^{-\frac{\bc}{d} \cdot \bx}\Delta u_i + \bc \cdot  e^{-\frac{\bc}{d} \cdot \bx}\nabla u_i + \frac{e^{-\frac{\bc}{d} \cdot \bx}}{\tau}u_i = \tilde{g}_i,
\end{equation}
where, $\tilde{g}_i = G( e^{-\frac{\bc}{d} \cdot \bx}\tilde{u}_i) + \frac{ e^{-\frac{\bc}{d} \cdot \bx}}{\tau}\tilde{u}_i.$ We remind the reader that formally, we have that $w_i^{n} = e^{-\frac{\bc}{d} \cdot \bx}\tilde{u}_i,$ and $w_i^{n+1} = e^{-\frac{\bc}{d} \cdot \bx}u_i.$ We note that 
$$\nabla \cdot \left(-d  e^{-\frac{\bc}{d} \cdot \bx} \nabla u_i\right) = -d  e^{-\frac{\bc}{d} \cdot \bx}\Delta u_i + \bc \cdot  e^{-\frac{\bc}{d} \cdot \bx}\nabla u_i.$$
Thus, we can write Equation (\ref{eq:notDivForm}) as 
\begin{equation}\label{eq:divForm}
\nabla \cdot \left(-d  e^{-\frac{\bc}{d} \cdot \bx} \nabla u_i\right) + \frac{e^{-\frac{\bc}{d} \cdot \bx}}{\tau}u_i = \tilde{g}_i.
\end{equation}
The jump interface conditions (\ref{eq:densityConditionSemiDiscreteTrunc}) become 
\begin{equation}\label{eq:divFormJump}
e^{-\frac{\bc}{d} \cdot \bx}u_0 =  \kappa e^{-\frac{\bc}{d} \cdot \bx} u_1 \implies u_0 =  \kappa u_1.
\end{equation}
The flux interface conditions (\ref{eq:fluxConditionSemiDiscreteTrunc}) become
\begin{equation}\label{eq:divFormFlux}
de^{-\frac{\bc}{d} \cdot \bx}\partial_{\bnu}u_0 = de^{-\frac{\bc}{d} \cdot \bx}\partial_{\bnu}u_1 \implies \partial_{\bnu}u_0=  \partial_{\bnu}u_1.
\end{equation}
\begin{theorem}
When $\bc \neq 0$ and $d_1 = d_0 = d,$ then Problem ($\text{\textbf{P}}'$) is well-posed. Moreover, $w \in \mathcal{V}_{\kappa}$ is the unique solution of System (\ref{eq:2DSemiDiscreteTrunc}) if and only if it is the unique solution of Problem ($\text{\textbf{P}}'$). 
\end{theorem}
\begin{proof}
To prove this statement, we look at the following variational problem: Find $u \in \mathcal{V}_{\kappa}$ such that 
$$\tilde{a}_d(u, v) = \tilde{\ell}_d(v), \qquad \forall v \in \mathcal{V}_{\kappa},$$
where
\begin{equation*}
\tilde{a}_d(u, v) = \int_{\Omega_0} de^{-\frac{\bc}{d} \cdot \bx} \nabla u_0 \nabla v_0 + \frac{e^{-\frac{\bc}{d} \cdot \bx}}{\tau}u_0v_0 \ \rd \bx 
+ \int_{\Omega_1} \kappa\left(de^{-\frac{\bc}{d} \cdot \bx} \nabla u_1 \nabla v_1 + \frac{e^{-\frac{\bc}{d} \cdot \bx}}{\tau}u_1v_1\right) \ \rd \bx,
\end{equation*}
and
$$
\tilde{\ell}_d(v) =  \int_{\Omega_0} \tilde{g}_0 v_0 \ \rd\bx + \int_{\Omega_1} \kappa \tilde{g}_1 v_1 \ \rd\bx. 
$$
Since $e^{-\frac{\bc}{d} \cdot \bx}$ is bounded and positive on $\Omega_i,$ $i = 0, 1,$ this bilinear form is coercive and bounded. Thus, by the Lax--Milgram Lemma, this variational problem is well-posed. Then, by following similar steps as in the proof of Theorem \ref{theo:c0}, we find that solutions to this variational problem are solutions to System (\ref{eq:2DSemiDiscreteTrunc}) and vice-versa. Thus, there is a unique solution to System (\ref{eq:2DSemiDiscreteTrunc}), when $\bc \neq 0$ and $d_1 = d_0.$

Thus, Problem ($\text{\textbf{P}}'$) has a unique solution when $\bc \neq 0$ and $d_1 = d_0$ and 
$$\left\|w\right\|_X \leq C\left\|g\right\|_{0,\Omega}, \ \ C > 0 \text{ is a constant.}$$
%\qed
\end{proof}
The above two theorems prove that Problem ($\text{\textbf{P}}'$) and the strong formulation are well-posed in the two cases: 1) $\bc = 0;$ 2) $\bc \neq 0$ and $d_1 = d_0.$ For the general case, we are unable to show that the bilinear form is coercive. One can try instead to prove that the weaker conditions of the Banach-Ne\u{c}as-Babu\u{s}ka Theorem \cite{Ern:2004:BookTheoryAndPractise} are satisfied. We leave this challenge for a future endeavour. For the rest of this chapter, we will assume that Problem ($\text{\textbf{P}}'$) is well-posed in the general case.  
\subsection{Existence and Uniqueness for Problem (\textbf{P})}
Finally, we show that Problem (\textbf{P}) is well-posed. 

\begin{theorem}\label{thm:PUnique}
Problem (\textbf{P}) has a unique solution.
\end{theorem}
\begin{proof}
We have that
\begin{align*}
w \text{ is a solution of Problem (\textbf{P})  } &\iff w \text{ is a solution of System (\ref{eq:2DSemiDiscreteTrunc}). }
\\
&\iff w \text{ is a solution of Problem ($\text{\textbf{P}}'$)}
\end{align*}
Since Problem ($\text{\textbf{P}}'$) has a unique solution, we conclude that Problem (\textbf{P}) has a unique solution.
%\qed
\end{proof}

We now have existence and uniqueness for Problem (\textbf{P}). Thus, by the Banach-Ne\u{c}as-Babu\u{s}ka Theorem \cite{Ern:2004:BookTheoryAndPractise} the following conditions are satisfied:
\begin{align}
&\exists \sigma > 0, \qquad \inf_{w \in \mathcal{V}_\kappa}\sup_{v \in \mathcal{V}_1} \frac{a(w, v)}{\left\|w\right\|_{X}\left\|v\right\|_{X}} \geq \sigma,
\\
&\forall v\in\mathcal{V}_1, \ \ (\forall w \in \mathcal{V}_\kappa, a(w, v) = 0) \ \ \implies (v = 0).
\end{align}
Moreover, we have the following \textit{a priori} estimate:
\begin{equation}
\forall g \in \mathcal{V}_1', \ \ \left\|w\right\|_{\mathcal{V}_\kappa} \leq \frac{1}{\sigma}\left\|g\right\|_{\mathcal{V}_1'},
\end{equation}
where $\mathcal{V}_1'$ denotes the dual of $\mathcal{V}_1.$ In particular, this estimate holds for $g$ as defined in System (\ref{eq:2DSemiDiscreteTrunc}).
\subsection{Hybrid Variational Formulation}\label{sec:hybridVariationalForm}
In the previous section, we introduced Problem (\textbf{P}), a variational formulation where the jump condition along the interface $\Gamma$ is explicitly imposed in the functional space $\mathcal{V}_{\kappa}.$ With a hybrid formulation, the density jump across the interface is weakly imposed through a Lagrange multiplier, which belongs to the following space:
$$M = \{\mu \in  H^{-1/2}(\Gamma) : \exists \bq \in H(\diver;\Omega) \text{ such that }  \bq\cdot\bm{\nu} = \mu \text{ on } \Gamma\},$$
where $ H(\diver;\Omega) = \{\bq \in \left(L^2(\Omega)\right)^2 :  \diver \bq \in L^2(\Omega)\}$ \cite{Raviart:1977:MathComput}. With this functional framework, we define a hybrid formulation, which we call Problem (\textbf{Q}):
\begin{description}
\item[(\textbf{Q})]\label{prob:hybrid}
Find $(w, \lambda) \in X\times M$ such that 
\begin{align}\label{eq:hybridFirstAbstractDirichlet}
&a(w, v) + b_1(v, \lambda) = \langle g,v \rangle_{X', X}, & \forall v \in X, 
\\ \label{eq:hybridSecondAbstractDirichlet}
&b_{\kappa}(w, \mu) = 0, & \forall \mu \in M, 
\end{align} 
\end{description}
where $a(w, v)$ is defined as in Equation (\ref{eq:a(w,v)}) and 
\begin{equation*}
b_{\kappa}(w, \lambda) = \int_\Gamma \lambda(w_0 -  \kappa w_1) \ \rd s,
\end{equation*}
thus, 
\begin{equation*}
b_1(w, \lambda) = \int_\Gamma \lambda(w_0 -  w_1)\ \rd s.
\end{equation*}
Here, the integral $\int_{\Gamma}$ represents the duality product between $H^{-1/2}(\Gamma)$ and $H^{1/2}(\Gamma);$ i.e., $\int_\Gamma \lambda(w_0 -  \kappa w_1)\ \rd s =  \langle \lambda, w_0 - \kappa w_1\rangle_{H^{-1/2}, H^{1/2}},$ and $w_i$ denotes $w_i|_{\Gamma} = \gamma (w_i)$ where $\gamma: H^1(\Omega_i) \to H^{1/2}(\Gamma)$ is the trace operator.  

Under the assumption that Problem (\textbf{P}) is well-posed in the most general case, that is for nonzero shifts and differing diffusion rates, we now prove that Problem (\textbf{Q}) is well-posed in this same general setting.
\subsection{Existence and Uniqueness of Solutions for the Hybrid Variational Form}\label{sec:existenceAndUniqueness}
Here we prove the existence and uniqueness of solutions to Problem (\textbf{Q}). Before we present our proof, we establish characterisations of our variational space $\mathcal{V}_{\kappa}$ parameterised by $\kappa > 0$ through $b_\kappa.$

\begin{lemma}\label{lem:uniqueMu2}
A continuous linear functional $L$ on the space $X$ vanishes on $\mathcal{V}_{\kappa}$ if and only if there exists a unique element $\mu \in M$ such that %
$$\forall v \in X, \ L(v) = \int_\Gamma \mu(v_0 - \kappa v_1) \ \rd s.$$
\end{lemma}
In the case that $\kappa = 1,$ the proof of Lemma \ref{lem:uniqueMu2} is identical to the proof of Lemma 1 in \cite{Raviart:1977:MathComput}. The proof in the case where $\kappa \neq 1$ requires more than minor modifications, so we present it here.
%%
%~\ \\ 
%~ \\
%%
\begin{proof}
By the Riesz Representation Theorem and the Hahn-Banach Theorem, there exist $q_j \in L^2(\Omega_i)$, $j = 0,1,2,$ such that any continuous linear function on $H^1(\Omega_i)$ can be written in the following way 
$$v \mapsto \int_{\Omega_i}\kappa_i\left(q_1\frac{\partial v}{\partial x} + q_2\frac{\partial v}{\partial y} + q_0v\right)\ \rd\bx,$$
where 
$$\kappa_i = \begin{cases} 1, \ \  \bx \in \Omega_0,
\\ 
\kappa, \ \ \bx \in \Omega_1.
\end{cases}
$$ 

Hence, given a continuous linear functional $L$ on $X$, there exist $q_j \in L^2(\Omega)$, $j = 0,1,2$ such that 
$$\forall v \in X, \ L(v) = \sum_{i = 0}^1\int_{\Omega_i}\kappa_i\left(q_1\frac{\partial v}{\partial x} + q_2\frac{\partial v}{\partial y} + q_0v\right)\ \rd\bx.$$ 

Assume that $L(v)$ vanishes on $\mathcal{V}_{\kappa};$ i.e., 
\begin{equation}\label{eq:vanish}
\forall v \in \mathcal{V}_{\kappa}, \ L(v) = \sum_{i = 0}^1\int_{\Omega_i}\kappa_i\left(q_1\frac{\partial v}{\partial x} + q_2\frac{\partial v}{\partial y} + q_0v\right)\ \rd\bx = 0.
\end{equation}
Then, for all $v \in \mathcal{D}(\Omega_i) \subset \mathcal{V}_{\kappa}$ we have that  
\begin{equation} 
\int_{\Omega_i}\kappa_i\left(q_1\frac{\partial v}{\partial x} + q_2\frac{\partial v}{\partial y}\right)\ \rd\bx  = -\int_{\Omega_i}\kappa_iq_0v\ \rd\bx.
\end{equation}
Using the definition of derivatives in the sense of distributions, we obtain
\begin{equation*} 
\left\langle\kappa_i\left(\nabla \cdot \bq\right), v\right\rangle  = \int_{\Omega_i}\kappa_iq_0v\ \rd\bx.
\end{equation*}
Thus, inside each $\Omega_i,$ $ \diver \bq = q_0 \in L^2(\Omega_i)$ in the sense of distributions, where $\bq = (q_1, q_2).$ Therefore, $\bq \in H(\diver; \Omega_i),$ and the following Green's formula holds within each $\Omega_i$:
\begin{equation}\label{eq:greensHDiv}
\forall v \in H^1(\Omega_i), \ \ \int_{\Omega_i} \nabla v_i \cdot \bq+ v_i \ \diver\bq\ \rd\bx = \int_{\Gamma} v_i\bq \cdot \bnu_i \ \rd s.
\end{equation}
We replace $q_0$ with $\diver \bq$ in Equation (\ref{eq:vanish}). Let $\gamma^*(\bq|_{\Omega_i}) = \bq_i\cdot\bm{\nu_i},$ where $\gamma^*: H(\diver;\Omega_i) \to H^{-1/2}(\Gamma)$ and $\bm{\nu_i}$ is the outward pointing normal to the boundary of $\Omega_i.$ Then for all $v \in \mathcal{V}_{\kappa}$, by applying Green's Formula (\ref{eq:greensHDiv}) within each $\Omega_i$, we find that 
\begin{equation*}
\sum_{i = 0}^1\int_{\Omega_i}\kappa_i\left( \nabla v \cdot \bq+ v \ \diver\bq\right)\ \rd\bx = 0 \iff \int_\Gamma \kappa v_1\left(\bq_1\cdot\bnu_1 + \bq_0\cdot\bnu_0\right)  \ \rd s = 0.
\end{equation*}
This equivalence follows since on $\Gamma,$ $v_0 = \kappa v_1.$ Since this relationship is true for all $v_1 \in H^{1/2}(\Gamma),$ we have that in $H^{-1/2}(\Gamma)$
\begin{equation*}
\bq_1\cdot\bnu_1 + \bq_0\cdot\bnu_0= 0 \iff \bq_0\cdot\bnu = \bq_1\cdot\bnu,
\end{equation*}
where equivalence holds as $\nu := \nu_0 = -\nu_1.$ Thus, $\bq$ belongs to $H(\diver; \Omega).$ Indeed, take $\phi \in \mathcal{D}(\Omega).$ Then the following holds:
\begin{align*}
\int_{\Omega} \bq \nabla \phi \ \rd\bx &= \sum_{i = 0}^1\left[-\int_{\Omega_i} \phi \diver\bq \ \rd \bx + \int_{\Gamma}\phi \bq \cdot \bnu_i \ \rd s\right] = -\int_{\Omega} \phi \diver\bq \ \rd \bx.
\end{align*}
Therefore, $\diver \bq \in L^2(\Omega)$ and we can write 
\begin{equation}\label{eq:diverRep}
\forall v \in X, \ L(v) = \sum_{i = 0}^1\int_{\Omega_i} \kappa_i\left(\nabla v \cdot \bq+ v \ \diver \bq\right)\ \rd\bx.
\end{equation}
Conversely, any linear functional of the form (\ref{eq:diverRep}) is continuous on $X$ and vanishes on $\mathcal{V}_{\kappa}.$

Using Green's formula on Equation (\ref{eq:diverRep}) and the definition of $H^{-1/2}(\Omega_i)$ within each $\Omega_i$ we obtain for all $v \in X$
$$
L(v) = \sum_{i = 0}^1\int_\Gamma \kappa_i v_i\bq\cdot\bnu_i\ \rd s = \int_\Gamma\bq\cdot\bnu\left(v_0 - \kappa v_1\right)\ \rd s =  \int_\Gamma\mu\left(v_0 - \kappa v_1\right)\ \rd s, \ \ \mu \in M.
$$
The function $\bq$ is not uniquely determined, but the corresponding element $\mu \in M$ is unique. Indeed, assume that 
$$\forall v \in X, \ \int_\Gamma\mu\left(v_0 - \kappa v_1\right)\ \rd s = 0.$$ 
Then, by taking $v_{1} \equiv 0$ on $\Omega_1$ and $v_{0} \in H^1(\Omega_0)$ such that $v_0 \not\equiv 0$ in $H^{1/2}(\Gamma),$ we find that 
$$\int_{\Gamma}\mu v_0 \ \rd s = 0.$$
Thus, we have that
$$\forall v \in H^1(\Omega_0), \quad \int_{\Gamma}\mu v\ \rd s = 0,$$ 
which implies that $\mu = 0$ on $\Gamma$ by the surjectivity of the trace operator $v \mapsto v|_{\Gamma}$ from $H^1(\Omega_0)$ onto $H^{1/2}(\Gamma).$
%\qed
\end{proof}
As a consequence of Lemma \ref{lem:uniqueMu2}, we can characterise $\mathcal{V}_1$ and $\mathcal{V}_{\kappa}$ as subspaces of $X,$ in the following ways:
\begin{equation}\label{eq:characterV1}
\mathcal{V}_1 = \{ v \in X \ : \ \forall \mu\in M, \ b_1(v, \mu) = 0\},
\end{equation}
and
\begin{equation}\label{eq:characterV2}
\mathcal{V}_{\kappa} = \{ v \in X \ : \ \forall \mu\in M, \ b_{\kappa}(v, \mu) = 0\}.
\end{equation}
\begin{theorem} 
Under the assumption that Problem (\textbf{P}) is well-posed in the general case (that is when $\bc\neq0$ and $d_0 \neq d_1$) Problem (\textbf{Q}) has a unique solution $(w, \lambda) \in X \times M.$ Moreover, $w \in \mathcal{V}_{\kappa}$ is the solution of System (\ref{eq:2DSemiDiscreteTrunc}) and we have that
\begin{equation*}
\lambda = -\left(d_i\partial_{\bnu}w_i + (\bc\cdot \bnu)w_i\right) \ \text { on } \ \Gamma \ \text{ for } \ i = 0, 1. 
\end{equation*}
\end{theorem}
\begin{proof}
Let $(w, \lambda) \in X \times M$ be a solution to Problem (\textbf{Q}). Then from the characterisation of $\mathcal{V}_{\kappa},$ Equation (\ref{eq:characterV2}), we see that $w\in\mathcal{V}_{\kappa}.$ Choosing $v\in\mathcal{V}_1$, from the characterisation of $\mathcal{V}_1,$ Equation (\ref{eq:characterV1}), we obtain
$$a(w,v) = \int_{\Omega} gv \ \rd \bx \bm, \ \forall v \in \mathcal{V}_1.$$
This equation is nothing but Problem (\textbf{P}). Thus, existence and uniqueness of $w \in \mathcal{V}_{\kappa}$ follows from the well-posedness of Problem (\textbf{P}). Additionally, $w \in \mathcal{V}_{\kappa}$ is then a solution of System (\ref{eq:2DSemiDiscreteTrunc}). 

Now, suppose that $w \in \mathcal{V}_{\kappa}$ is a solution of System (\ref{eq:2DSemiDiscreteTrunc}). Consider the linear continuous function on $X$
$$L(v) = \int_\Omega gv \ \rd \bx - a(w, v).$$
The right-hand side is Problem (\textbf{P}). Hence, for $v \in \mathcal{V}_1,$ $L(v)$ vanishes. Thus, by Lemma \ref{lem:uniqueMu2}, there exists a unique $\lambda \in M$ such that 
\begin{equation*}
\forall v \in X, \ b_1(v, \lambda) = \int_{\Omega} gv \ \rd\bx - a(w,v).
\end{equation*}
Additionally, since $w\in \mathcal{V}_{\kappa},$ we have that $b_{\kappa}(w, \lambda) = 0.$ Thus, $(w, \lambda)$ is the solution to the hybrid formulation. 

Finally, since $g_i = \frac{w_i}{\tau} - d_i\Delta w_i - (\bc\cdot \nabla)w_i$ within each $\Omega_i$, by applying Green's formula in each $\Omega_i$ with $\nabla w = \bq,$ we find that, $\forall v \in X$
\begin{align*}
b_1(v, \lambda) &= \sum_{i = 0}^1\int_{\Omega_i}\left(-d_i\Delta w_i  - (\bc\cdot \nabla)w_i + \frac{1}{\tau}w_i\right)v_i \ \rd\bx - a(u, v) 
\\
&= \int_{\Gamma} \left(d_1\partial_{\bnu}w_1 + (\bc\cdot \bnu)w_1\right)v_1 - \left(d_0\partial_{\bnu}w_0 + (\bc\cdot \bnu)w_0 \right) v_0 \ \rd s,
\\
&= \int_{\Gamma}(v_0 - v_1)\lambda\ \rd s.
\end{align*}
Thus, we have that on $\Gamma,$ 
$$\lambda = -\left(d_1\partial_{\bnu}w_1 + (\bc\cdot \bnu)w_1\right) = -\left(d_0\partial_{\bnu}w_0 + (\bc\cdot \bnu)w_0\right).$$
%\qed
\end{proof}
In the present section, we introduced a hybrid variational formulation, called Problem (\textbf{Q}). The hybrid system has two unknowns, a primal variable, $w,$ and a dual variable, $\lambda.$ We proved that Problem (\textbf{Q}) has a unique solution. Furthermore, we proved that the primal variable is the unique solution for System (\ref{eq:2DSemiDiscreteTrunc}), and that the dual variable is the flux across the interface $\Gamma.$ We are now ready to construct a finite element method to approximate the solution of Problem (\textbf{Q}). 
\section{Finite Element Formulation}\label{sec:finiteElementFormulation}
Our goal now is to construct a finite element method as a solver for Problem (\textbf{Q}). In this section, we assume that $\Gamma$ is a polygon and denote by $\Gamma_j,$ $j = 1,...,J,$ the line segments that form $\Gamma.$ This assumption is necessary for the construction of our proofs, yet in implementation, we do not consider it a necessary assumption to obtain accurate results. 

We denote by $\mathbb{P}_s$ the space of polynomials of degree $\leq s.$ Let $\mathcal{T}_h^i$ denote a triangulation of $\Omega_i,$ where any element $K$ of $\mathcal{T}_h^i$ is triangular. Here, $h$ denotes the length of the longest edge of all the elements $K.$ Then, we can define the finite-dimensional spaces
$$\mathcal{U}_h^i = \mathcal{U}_h(\Omega_i) =  \{v_h \in \mathcal{C}(\Omega_i) : v_h|_{K} \in \mathbb{P}_s(K), \ \forall K \in \mathcal{T}_h^i, \ v_h|_{\partial\Omega_i\cap\partial\Omega} = 0\}.$$
We consider the subset of $X:$
$$X_h = \mathcal{U}_h^0 \times \mathcal{U}_h^1.$$

The trace of the triangulation within each $\Omega_i$ over $\Gamma$ creates a triangulation, $t_i,$ $i = 0, 1,$ along $\Gamma.$ We denote by $F$ the elements of $t_i.$ The distinction between these two triangulations is relevant since it is not necessary for the triangulations to conform across $\Gamma.$ With this notation, we define the finite-dimensional space
$$M_h = \{\mu_h \in \mathcal{C}(\Gamma) : \forall F \in t_0, \  \mu_h|_{F} \in \mathbb{P}_s(F)\}.$$
Here, we do not make any reduction on the polynomial order near vertices of polygonal domains, as has been done in comparable decomposition methods; see for example \cite{Bernardi:1994:ANewConform}.
\begin{remark}
We choose, in line with our numerical implementation, to define $M_h$ through the triangulation resulting from $\mathcal{T}_h^0.$ It is possible to define $M_h$ in other ways, for example, by the triangulation resulting from $\mathcal{T}_h^1,$ or by the intersection of the two \cite{Belgacem:1999:NumerMath}.
\end{remark}
Our finite-dimensional approximation of the hybrid formulation, which we call Problem ($\text{\textbf{Q}}_h$), is: 
\begin{description}
\item[($\text{\textbf{Q}}_h$)]\label{prob:hybridApprox}
Find $(w_h, \lambda_h) \in X_h\times M_h$ such that 
\begin{align}\label{eq:hybridApproxA}
&a(w_h, v_h) + b_1(v_h, \lambda_h) = \langle g,v_h\rangle, & \forall v_h \in X_h, 
\\ \label{eq:hybridApproxB}
&b_{\kappa}(w_h, \mu_h) = 0, & \forall \mu_h \in M_h.
\end{align}
\end{description}
For well-posedness of our finite element method, we need the following spaces:
$$\mathcal{V}_{\kappa h} = \{v_h \in X_h \ : \ b_{\kappa}(v_h, \mu_h) = 0, \  \forall \mu_h \in M_h\},$$
and
$$\mathcal{V}_{1h} = \{v_h \in X_h \ : \ b_{1}(v_h, \mu_h) = 0, \  \forall \mu_h \in M_h\};$$
and auxiliary problems:
\begin{description}
\item[($\text{\textbf{P}}_h$)]\label{prob:Ph}
Find $w_h \in \mathcal{V}_{\kappa h}$ such that 
\begin{align}\label{eq:Ph}
&a(w_h, v_h) = \ell(v_h), & \forall v_h \in \mathcal{V}_{1h},
\end{align}
\end{description}
and
\begin{description}
\item[($\text{\textbf{P}}_h'$)]\label{prob:Ph'}
Find $w_h \in \mathcal{V}_{\kappa h}$ such that 
\begin{align}\label{eq:Ph'}
&\tilde{a}(w_h, v_h) = \tilde{\ell}(v_h), & \forall v_h \in \mathcal{V}_{\kappa h}.
\end{align}
\end{description}
\begin{lemma}\label{lem:PhequivPh'}
The function $w_h \in \mathcal{V}_{\kappa h}$ is a solution to Problem ($\text{\textbf{P}}_h$) if and only if it is a solution to Problem ($\text{\textbf{P}}_h'$). 
\end{lemma}
\begin{proof}
Suppose that $w_h \in \mathcal{V}_{\kappa h}$ is a solution to Problem ($\text{\textbf{P}}_h$). Let $(v_{0h}, v_{1h}) \in \mathcal{V}_{1h}.$ Then, the equation 
$$v_h = (v_{0h}, v_{1h}) = (u_{0h}, \kappa u_{1h})$$
defines a one-to-one relation between $\mathcal{V}_{1h}$ and  $\mathcal{V}_{\kappa h}.$ Upon substituting $u_h$ into Equation (\ref{eq:Ph}) we see that $w_h$ satisfies Equation (\ref{eq:Ph'}) for any $v_h$ defined through this relation. 

Arguing in a similar manner in the other direction proves that $w_h,$ a solution of Problem ($\text{\textbf{P}}_h'$), is also a solution of
Problem ($\text{\textbf{P}}_h$), thus concluding the proof. 
%\qed
\end{proof}
For the rest of this section, we assume that $\bc = 0.$ We use the 0-subscript notation, as done in Section \ref{sec:existAndUniqeProbP'}, to indicate that this assumption holds. In later sections, we validate our scheme for cases where $\bc \neq 0.$ An analytical proof of existence and uniqueness for the discretised system in the most general case remains an open question. 

Before proving well-posedness of our finite-dimensional problems, we prove that the compatibility condition (3.7) in \cite{Raviart:1977:MathComput} holds. 
\begin{proposition}\label{prop:setZero}
For any $j = \kappa, 1,$ $\{\mu_h \in M_h \ ; \ \forall v_h \in X_h, \ b_{j}(v_h, \mu_h) = 0\} = \{0\}.$
\end{proposition}
\begin{proof}
Consider $\mu_h \in M_h$ and choose $v_h \in X_h$ such that $v_{1h}|_{\Gamma} = 0$ and $v_{0h}|_{\Gamma} = \mu_h.$ Then 
$$b_{j}(v_h, \mu_h) = \int_{\Gamma}\mu_h^2 \ \rd s = 0 \iff \mu_h = 0.$$
%\qed
\end{proof}
\begin{theorem}\label{theo:uniqueSolutionsApprox}
\begin{enumerate}
\item Problem ($\text{\textbf{P}}_{h0}'$) has a unique solution. 
\item Problem ($\text{\textbf{P}}_{h0}$) has a unique solution. 
\item Problem ($\text{\textbf{Q}}_{h0}$) has a unique solution. 
\end{enumerate}
\end{theorem}
\begin{proof}
The proof is similar to that of Theorem 2 of \cite{Raviart:1977:MathComput}.
\begin{enumerate}
\item We have that $V_{\kappa h} \subset X_h \subset X.$ Since $\tilde{a}_0(v_h, v_h)$ is coercive and bounded over $X,$ we conclude that $\tilde{a}_0(v_h, v_h)$ is coercive and bounded over $\mathcal{V}_{\kappa h}.$ Thus, by the Lax--Milgram Lemma, Problem ($\text{\textbf{P}}_{h0}'$) has a unique solution.
\item By Lemma \ref{lem:PhequivPh'}, Problem ($\text{\textbf{P}}_{h0}$) has a unique solution.
\item Problem ($\text{\textbf{Q}}_{h0}$) is equivalent to an $N \times N$ linear system with $N = \dim(X_h) + \dim(M_h).$ Thus, existence of a solution follows from uniqueness, so that we need only to consider the corresponding homogeneous equation. We assume that $g = 0.$ Then, necessarily, from Problem ($\text{\textbf{P}}_{h0}$), $w_h = 0.$ We have then that $\lambda_h$ satisfies
$$b_1(v_h, \lambda_h) = 0, \qquad \forall v_h \in X_h.$$
Hence, from Proposition \ref{prop:setZero}, $\lambda_h = 0.$
\end{enumerate}
%\qed
\end{proof}
With uniqueness and existence proven, we note that in fact, if $(w_h, \lambda_h) \in X_h \times M_h$ is the unique solution of Problem ($\text{\textbf{Q}}_{h0}$), then $w_h \in \mathcal{V}_{\kappa h}$ is the unique solution of Problem ($\text{\textbf{P}}_{h0}$). 
\section{Error Analysis}\label{sec:errorAnalysis}
In this section, we derive error estimates on $w - w_h.$ Particularly novel are the optimal convergence rates of the primal variable given nonstandard jump conditions. Moreover, contrary to the existing literature, is the absence of restrictions on the discretisation of the dual variable $\lambda$ near the corners of the interface between subdomains. Throughout this section, we assume $\bc = 0.$ However, in Section \ref{sec:2dNumerics}, we relax this assumption and demonstrate that even for nonzero shifts, our results remain relevant. Our overarching goal is to approximate well the solution $w$ of System (\ref{eq:2DSemiDiscreteTrunc}). We leave the error analysis of $\lambda$ for future work. 

We introduce the auxilliary hybrid finite element formulation:
\begin{description}
\item[($\text{\textbf{Q}}_{h0}'$)]\label{prob:hybridApproxPrime}
Find $(w_h, \lambda_h) \in X_h\times M_h$ such that 
\begin{align}\label{eq:hybridApproxPrime}
&\tilde{a}_0(w_h, v_h) + b_\kappa(v_h, \lambda_h) = \tilde{\ell}(v_h), & \forall v_h \in X_h, 
\\ \label{eq:hybridApproxPrime}
&b_{\kappa}(w_h, \mu_h) = 0, & \forall \mu_h \in M_h.
\end{align}
\end{description}

In the same way that Problem ($\text{\textbf{Q}}_{h0}$) has a unique solution, we can show that Problem ($\text{\textbf{Q}}_{h0}'$) has a unique solution. Additionally, whenever $(w_h, \lambda_h) \in X_h\times M_h$ is the unique solution to Problem ($\text{\textbf{Q}}_{h0}'$), then $w_h \in \mathcal{V}_{\kappa h}$ is also the unique solution to Problem ($\text{\textbf{P}}_{h0}'$). 

We define the following norm: 
$$ \left\|v_h\right\|_a = (\tilde{a}_0(v_h, v_h))^{1/2},$$
for which there exist $C_*$ and $C^*$ dependent only on $\kappa,$ $d_i$ and $\tau$ such that for $v_h \in X_h$
$$
C_*\left\| v_h \right\|_a \leq \left\| v_h \right\|_X \leq C^*\left\| v_h \right\|_a.
$$
We denote by $\left\| \cdot \right\|_{a, \Omega_i}$ the norm 
$$
\left\| w \right\|_{a, \Omega_i} = \left(\int_{\Omega_i}d_i |\nabla w_i| ^2 \ \rd \bx + \int_{\Omega_i} \frac{1}{\tau}|w_i| ^2 \ \rd \bx\right)^{\frac{1}{2}}. 
$$
\begin{theorem}
The solution $w_h \in \mathcal{V}_{\kappa h}$ of Problem ($\text{\textbf{P}}_{h0}'$) satisfies
\begin{equation}\label{eq:approxAndConsist}
\left\|w - w_h\right\|_a = \left( \left(\inf_{v_h \in \mathcal{V}_{\kappa h}} \left\|w - v_h\right\|_a\right)^2 + \left(\inf_{\mu_h \in M_h} \sup_{v_h \in \mathcal{V}_{\kappa h}} \frac{b_{\kappa}(v_h, \lambda - \mu_h)}{\left\|v_h\right\|_a} \right)^2 \right)^{1/2}.
\end{equation}
\end{theorem}
\begin{proof}
The result is a direct consequence of Theorem 3 of \cite{Raviart:1977:MathComput}.
%\qed
\end{proof}
The first term in Equation (\ref{eq:approxAndConsist}) is called the approximation error, the second term is called the consistency error. The consistency error results from $\mathcal{V}_{\kappa h} \not \subset \mathcal{V_\kappa};$ i.e., the finite-dimensional problem results in a nonconformal method. 

Before obtaining bounds on $\left\|w - w_h\right\|_a,$ we introduce some notation. We assume that $\mathcal{T}_h^i$ are regular families of triangulations in the sense that there exists a constant $\sigma_1 > 0$ independent of $h$ such that 
$$\max_{K \in \mathcal{T}_h^i} \frac{h_K}{\rho_K} \leq \sigma_1,$$
where $$\rho_K = \text{ diameter of the inscribed sphere in $K$.}$$

Additionally, we also assume that there exists a constant $\sigma_2$ such that for any two elements $F, F' \in t_0$ 
$$
|F| \geq \sigma_2 |F'|,
$$
where $|\cdot|$ denotes the length of the segment. 

\begin{theorem}\label{theo:consistError}
Let $M_h$ be associated to a regular family of triangulations. Let $\phi = (\phi_0, \phi_1) \in H^{\sigma+\frac{3}{2}}(\Omega_0)\times H^{\sigma+\frac{3}{2}}(\Omega_1),$ where $\frac{1}{2} \leq \sigma \leq s + 1$ such that $\frac{\partial \phi_0}{\partial \bnu} = \lambda \in M.$ Then there exists a constant $C,$ independent of $h$ such that 
\begin{equation}\label{eq:consistencyBound}
\inf_{\mu_h \in M_h} \sup_{v_h \in \mathcal{V}_{\kappa h}} \frac{b_{\kappa}(v_h, \lambda - \mu_h)}{\left\|v_h\right\|_a} \leq Ch^{\frac{1}{2} + \sigma}\left(\left\|\phi_0\right\|_{\sigma + \frac{3}{2}, \Omega_0} + \left\| \phi_1\right\|_{\sigma + \frac{3}{2}, \Omega_1}\right).
\end{equation}
\end{theorem}
\begin{proof}
We have that
\begin{align*}
b_{\kappa}(v_h, \lambda - \mu_h) &= \int_{\Gamma}(\lambda - \mu_h)(v_{0h} - \kappa v_{1h})\ \rd s
\\
&\leq \left\|\lambda - \mu_h\right\|_{H^{-1/2}(\Gamma)} \left(\left\|v_{0h}\right\|_{H^{1/2}(\Gamma)} + \left\| \kappa v_{1h} \right\|_{H^{1/2}(\Gamma)}\right)
\\
&\leq C \left\|\lambda - \mu_h\right\|_{H^{-1/2}(\Gamma)} \left(\left\|v_{0h}\right\|_{1, \Omega_0} + \left\| \kappa v_{1h} \right\|_{1, \Omega_1}\right)
\\
&\leq C \left\|\lambda - \mu_h\right\|_{H^{-1/2}(\Gamma)} \left\|v_{h}\right\|_a.
\end{align*}
Therefore, 
$$\frac{b_{\kappa}(v_h, \lambda - \mu_h)}{\left\|v_{h}\right\|_a} \leq C \left\|\lambda - \mu_h\right\|_{H^{-1/2}(\Gamma)},$$
which implies
$$\sup_{v_h\in\mathcal{V}_{\kappa h}}\frac{b_{\kappa}(v_h, \lambda - \mu_h)}{\left\|v_{h}\right\|_a} \leq C \left\|\lambda - \mu_h\right\|_{H^{-1/2}(\Gamma)}.$$
We can interpret $\mu_h$ as the orthogonal projection of $\lambda = \frac{\partial \phi}{\partial \bnu}$ from $L^{2}(\Gamma)$ onto $M_h.$ Thus, from Lemma 2.4 of \cite{Bernardi:1994:ANewConform}, we have that
\begin{equation*}
\left\|\lambda - \mu_h\right\|_{H^{-1/2}(\Gamma)}  \leq Ch_0^{\frac{1}{2} + \sigma}\left\| \frac{\partial \phi}{\partial \bnu}\right\|_{\sigma, \Gamma}
\leq Ch^{\frac{1}{2} + \sigma}\left(\left\|\phi_0\right\|_{\sigma + \frac{3}{2}, \Omega_0} + \left\|\phi_1\right\|_{\sigma + \frac{3}{2}, \Omega_1}\right). 
\end{equation*}
In the first line, $h_0$ is the longest edge of all the elements in $\Omega_0.$ The result then follows immediately. 
%\qed
\end{proof}
Before finding bounds on the approximation error, we construct a lifting from $M_h \to \mathcal{U}_h^0.$
\begin{lemma}\label{lem:lifting}
There exists a lifting $\mathcal{R}_{h, \Omega_0} : M_h \to \mathcal{U}_h^0,$ such that for any $\mu_h \in M_h,$ 
$$
\left\| \mathcal{R}_{h, \Omega_0}(\mu_h) \right\|_{a, \Omega_0} \leq C \left\| \mu_h \right\|_{\frac{1}{2}, \Gamma}. 
$$
\end{lemma}
\begin{proof}
Consider any $\mu_h \in M_h.$ Then, since $\mu_h$ is continuous along $\Gamma,$ $\mu_h \in H^1(\Gamma).$ Thus, by surjectivity of the trace operator, there exists an operator $R_{\Omega_0} : H^\frac{1}{2}(\Gamma) \to H^1(\Omega_0)$ such that 
$$
\left\| R_{\Omega_0}(\mu_h) \right\|_{1, \Omega_0} \leq C \left\| \mu_h \right\|_{\frac{1}{2}, \Gamma}. 
$$
We introduce the Scott-Zhang interpolator  $\SZ_h,$ (see expressions (2.12) and (2.13) in \cite{Scott:1990:MathComput}), where
$$
\SZ_h: H^1(\Omega_0) \to \mathcal{U}_h^0,
$$
and
$$
\left\| \SZ_h(v) \right\|_{1, \Omega_0} \leq C \left\| v \right\|_{1,\Omega_0}. 
$$
The Scott-Zhang interpolator preserves the boundary value; i.e., $\mathcal{SZ}_h(v)|_\Gamma = v|_\Gamma,$ if $v|_\Gamma \in M_h.$ 

We let $\mathcal{R}_{h, \Omega_0} = \SZ_h\circ R_{\Omega_0}.$ Then we have that, 
$$
\left\| \mathcal{R}_{h, \Omega_0} (\mu_h) \right\|_{a, \Omega_0} \leq C\left\| \mathcal{R}_{h, \Omega_0}(\mu_h) \right\|_{1, \Omega_0} \leq C\left\| R_{\Omega_0}(\mu_h) \right\|_{1, \Omega_0} \leq C \left\| \mu_h \right\|_{\frac{1}{2}, \Gamma}. 
$$
%
%\qed
\end{proof}
We now present bounds on the approximation error. 
\begin{theorem}\label{theo:approxError}
For $w = (w_0, w_1) \in H^{\delta}(\Omega_0) \times H^{\delta}(\Omega_1),$ with $1 < \delta \leq s+1,$ there exists a constant $C$ independent of $h$ such that 
\begin{equation}\label{eq:ApproxBound}
\inf_{v_h \in \mathcal{V}_{\kappa h}} \left\| w - v_h \right\|_a \leq C h^{\delta - 1}( \left\| w_0 \right\|_{\delta, \Omega_0} +  \left\| w_1 \right\|_{\delta, \Omega_1}).
\end{equation}
\end{theorem}

\begin{proof}
We have that 
\begin{align*}
\left\| w - v_h \right\|_a \leq \left\| w - v_h \right\|_{a, \Omega_0} + \sqrt{\kappa}\left\| w - v_h \right\|_{a, \Omega_1}.
\end{align*}
Let $\pi_{h, \Omega_i}$ be the Lagrange interpolant operator of degree $s$ on $\Omega_i.$ We can set $v_h|_{\Omega_1} = \pi_{h, \Omega_1}w_1, $ so that we have 
\begin{equation}\label{eq:errorInOmega1}
\left\| w -  \pi_{h, \Omega_1}w_1 \right\|_{a, \Omega_1}  \leq C\left\| w -  \pi_{h, \Omega_1}w_1 \right\|_{1, \Omega_1} \leq C h_1^{\delta - 1}\left\| w \right\|_{\delta, \Omega_1},
\end{equation}
where the last inequality comes from Theorem 3.1.6 in \cite{Ciarlet:2002:BookFiniteElement}. We need a suitable function $v_h= (\hat{w}_{0, h}, \hat{w}_{1, h})$ where $\hat{w}_{1, h} = \pi_{h, \Omega_1}w_1$ and $\hat{w}_{0, h}$ satisfies
\begin{equation}\label{eq:approxMustSatisfy}
\int_{\Gamma} \left(\hat{w}_{0, h} - \kappa\hat{w}_{1, h} \right)\mu_h \ \rd s = 0, \qquad \forall \mu_h \in M_h. 
\end{equation}
For any $\hat{\mu}_h \in M_h,$ we can define
$$
\hat{w}_{0, h} =  \pi_{h, \Omega_0}w_0 + \mathcal{R}_{h, \Omega_0}(\hat{\mu}_h).
$$
In order to satisfy Equation (\ref{eq:approxMustSatisfy}), we need
\begin{equation*}
\int_{\Gamma} \left(\hat{w}_{0, h} - \kappa\hat{w}_{1, h} \right)\mu_h \ \rd s = 0, \  \forall \mu_h \in M_h 
\iff \int_{\Gamma} \hat{\mu}_h\mu_h \ \rd s = \int_{\Gamma}\left(\kappa \pi_{h, \Omega_1}w_1 - \pi_{h, \Omega_0}w_0\right)\mu_h \ \rd s, \ \forall \mu_h \in M_h. 
\end{equation*}
Hence, we choose $\hat{\mu}_h$ to be the $L^2(\Gamma)$ projection of $\kappa \pi_{h, \Omega_1}w_1 - \pi_{h, \Omega_0}w_0$ on $M_h.$ Recalling that $b_{\kappa}(w, \mu_h) = 0,$ we find 
\begin{align}\label{eq:errorInOmega0}
\begin{split}
\left\| \hat{\mu}_h \right\|_{0, \Gamma} &\leq \left\| \kappa \pi_{h, \Omega_1}w_1 - \pi_{h, \Omega_0}w_0 + (w_0 - \kappa w_1) \right\|_{0, \Gamma}
\\
&\leq C\left[\left( \sum_{j = 1}^J \left\| \pi_{h, \Omega_1}w_1 - w_1 \right\|_{0, \Gamma_j}^2\right)^{1/2} +  \left( \sum_{j = 1}^J\left\| w_0 - \pi_{h, \Omega_0}w_0 \right\|_{0, \Gamma_j}^2\right)^{1/2}\right]
\\
&\leq C\left[\left( \sum_{j = 1}^J ch_1^{2\left(\delta - \frac{1}{2}\right)}|w_1|_{\delta - \frac{1}{2}, \Gamma_j}^2\right)^{1/2} +  \left( \sum_{j = 1}^Jch_0^{2\left(\delta - \frac{1}{2}\right)}|w_0|_{\delta - \frac{1}{2}, \Gamma_j}^2\right)^{1/2}\right]
\\
&\leq Ch^{\delta - \frac{1}{2}} \left[\left( \sum_{j = 1}^J \left\|w_1\right\|_{\delta - \frac{1}{2}, \Gamma_j}^2\right)^{1/2} +  \left( \sum_{j = 1}^J \left\|w_0\right\|_{\delta - \frac{1}{2}, \Gamma_j}^2\right)^{1/2}\right]
\\
&\leq Ch^{\delta - \frac{1}{2}} \left[\left(J \left\|w_1\right\|_{\delta, \Omega_1}^2\right)^{1/2} +  \left( J \left\|w_0\right\|_{\delta, \Omega_0}^2\right)^{1/2}\right] 
\\
&\leq Ch^{\delta - \frac{1}{2}} \left(\left\|w_1\right\|_{\delta, \Omega_1}+  \left\|w_0\right\|_{\delta, \Omega_0}\right).
\end{split}
\end{align}
In the third line, we applied Proposition 1.12 of \cite{Ern:2004:BookTheoryAndPractise}, while the second last inequality is due to the continuity of the trace operator (Theorem 1.5.2.1, \cite{Grisvard:2011:BookEllipticNonSmooth}). 

Thus, we have that 
\begin{align}\label{eq:infErrorBound}
\begin{split}
\left\| w - v_h \right\|_a &\leq \left\| w - \hat{w}_{0h} \right\|_{a, \Omega_0} + \sqrt{\kappa}\left\| w - \hat{w}_{1h} \right\|_{a, \Omega_1}
\\
&\leq C \left( \left\| w - \pi_{h,\Omega_0}w_0\right\|_{1, \Omega_0} + \left\| R_{h, \Omega_0}(\hat{\mu_h})\right\|_{1, \Omega_0} +  Ch_1^{\delta - 1}\left\| w_1 \right\|_{\delta, \Omega_1}\right)
\\
&\leq C \left(Ch_0^{\delta - 1}\left\| w_0 \right\|_{\delta, \Omega_0} + C\left\| \hat{\mu}_h\right\|_{\frac{1}{2}, \Gamma} +  Ch_1^{\delta - 1}\left\| w_1 \right\|_{\delta, \Omega_1}\right)
\\
&\leq C \left(Ch_0^{\delta - 1}\left\| w_0 \right\|_{\delta, \Omega_0} + Ch^{-\frac{1}{2}}\left\| \hat{\mu}_h\right\|_{0, \Gamma}  +  Ch_1^{\delta - 1}\left\| w_1 \right\|_{\delta, \Omega_1}\right)
\\
&\leq Ch^{\delta - 1} \left(\left\| w_0 \right\|_{\delta, \Omega_0} + \left\| w_1 \right\|_{\delta, \Omega_1}\right).
\end{split}
\end{align}
The second last line comes from the inverse inequality (Corollary 1.141, \cite{Ern:2004:BookTheoryAndPractise}). 

The result then follows immediately.  
%\qed
\end{proof}
We finally have arrived at the error bounds of $w - w_h.$
\begin{theorem}
For $(w_0, w_1) \in H^{\delta} \times H^{\delta}$ with $2 \leq \delta \leq s+1,$ there exists a constant $C$ independent of $h$ such that 
\begin{equation}
\left\| w - w_h \right\|_a \leq C h^{\delta - 1} \left(\left\| w_0 \right\|_{\delta, \Omega_0} + \left\| w_1 \right\|_{\delta, \Omega_1}\right).
\end{equation} 
\end{theorem}

\begin{proof}
The result is a direct consequence of Theorems \ref{theo:consistError} and \ref{theo:approxError}. 
%\qed
\end{proof}
\section{Numerical Test Cases}\label{sec:2dNumerics}
In the previous section, we provided theoretical \textit{a priori} error estimates in the case where $\bc = \bf{0}.$ In this section, we support our theoretical analysis with numerical test cases. We demonstrate that even with a non-zero shift ($\bc \neq 0$), we achieve optimal convergence rates in the $H^1$ norm and additionally second-order convergence rates in the $L^2$ norm. We then further validate our scheme with numerical test cases inspired by our previous study of a finite difference solver for a 1-dimensional moving-habitat model \cite{MacDonald:2021:MathBiosci}. Finally, from the perspective of application, we exemplify how the scheme can be used as a powerful tool to further the advancement of research in moving-habitat models. We implement our numerical solver using the FreeFEM$++$ finite element library \cite{Hecht:2012:JNumerMath}.
\subsection{Numerical Convergence}\label{sec:numericalConvergence}
We now check the order of convergence of our numerical method considering a non-zero shift; i.e., $\bc \neq \bf{0}.$ In our analysis, we consider two cases, the first being a conformal discretisation, the second being a nonconformal discretisation. In both cases, our computational domain is $\Omega_0 = \left[ 3,7\right]^2$ and $\Omega_1 = \left[-17, 19\right]\times\left[-17, 27\right]\backslash\bar{\Omega}_0.$ Solutions of System (\ref{eq:master2D}) with a constant shift approach a unique travelling pulse solution; this statement was proven analytically for globally continuous solutions \cite{Berestycki:2008:DiscretContinDynS}, and we suspect this to hold even with a jump in density across $\Gamma$ since the reaction term is monostable\footnote{In the sense that there is only one positive stable steady-state.}. In the reference frame, the travelling pulse translates to a steady-state solution; i.e., $\partial_t w = 0.$ 

\begin{figure}
\centering
\includegraphics[width = 6cm]{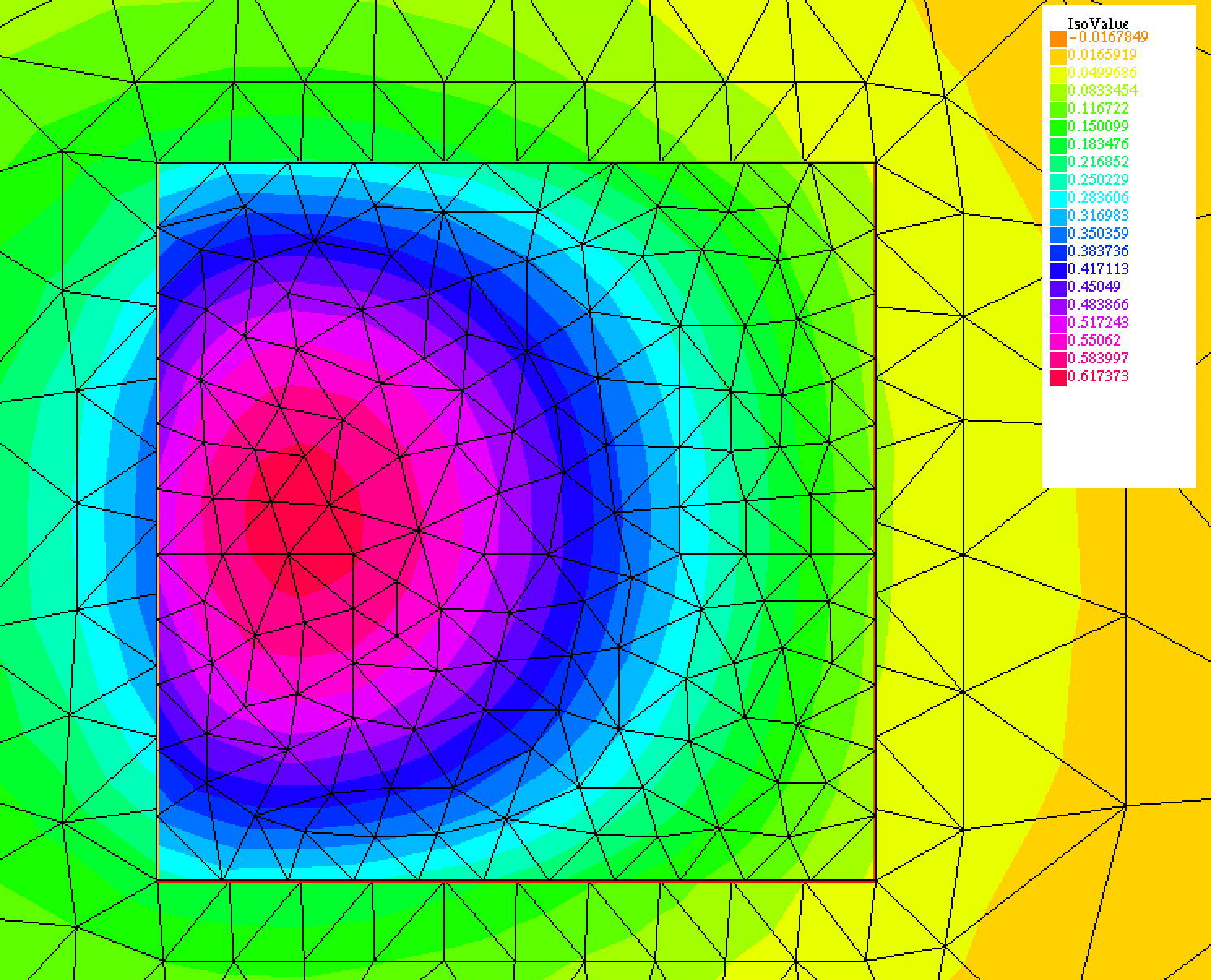}
\caption[ ]{The solution in and around the bounded domain over a nonconformal discretisation.}
\label{fig:nonConformalGrid}
\end{figure}

In lieu of an exact solution, we compare our numerical solutions to a reference solution computed over a conformal discretisation with 200 nodes on each side of the exterior rectangle and 400 nodes on each side of the interior square. We stop the time-stepping iterations when the solution has reached the travelling pulse profile, which we determine through the approximate time derivative in the $L^2$ norm. Our stopping condition is then:
\begin{equation}\label{eq:travellingPulseCond}
 \left\Vert\frac{w^{n+1} - w^{n}}{\tau}\right\Vert_{0, \Omega} = \left(\sum_{i = 0}^1 \int_{\Omega_i}\left(\frac{w^{n+1} - w^{n}}{\tau}\right)^2\right)^{\frac{1}{2}} < 10^{-5}.
\end{equation}
In our loop, we refine the time step, $\tau,$ so that Inequality (\ref{eq:travellingPulseCond}) holds for $\tau = 2.5\times10^{-2}.$ To calculate the integral in FreeFEM$++$, the software uses a $6^{th}$ order quadrature formula which is exact for polynomials of degree 5. The initial condition is the Gaussian function of height 10, with variance $(\sigma_x, \sigma_y) = (0.5, 0.5)$ and centred in $\Omega_0;$ i.e., 
\begin{equation}\label{eq:GaussInitial}
w^0(x,y) = \frac{10}{2\pi\sigma_x\sigma_y}\text{exp}\left(-\frac{1}{2}\left(\left(\frac{x - 5}{\sigma_x}\right)^2 + \left(\frac{y - 5}{\sigma_y}\right)^2\right)\right).
\end{equation}
We denote the number of nodes along $\Gamma$ from the exterior region, $\Omega_1,$ and the interior region, $\Omega_0,$ by $n_{\Gamma,1}$ and $n_{\Gamma, 0},$ respectively. In the conformal discretisation, $n_{\Gamma,0} = n_{\Gamma, 1}.$ In the nonconformal discretisation, $n_{\Gamma,0} + 1 = n_{\Gamma, 1}$  (Figure \ref{fig:nonConformalGrid}). We refine our grid by factors of two determined by $n_{\Gamma, 1}.$ We calculate the numerical solutions with a time step $\tau = 0.1.$ We stop the method when the approximate time derivative of the numerical solution is less than $10^{-5},$ as in Inequality (\ref{eq:travellingPulseCond}). The initial condition for each solution is the Gaussian function as in Equation (\ref{eq:GaussInitial}). We present results for two test cases. Test 1 has parameter values $\alpha = 0.5,$  $d_1 = 2,$  $d_0 = 1,$  $r = 1.2,$  $m = 1,$ $a = 0.8,$ and $ \bc  = (1,0).$ Test 2 has the same parameter values except with $\alpha = 0.7.$ These parameter values correspond to two qualitatively different solution shapes, where the first is hump shaped in the suitable habitat, while the second is a steeply decreasing function in the direction of the shift.
\begin{table}
\centering
\caption[]{Numerical Convergence in the $L^2$ norm and $H^1$ semi-norm}
\begin{tabular}{||c|c|| c || c | c ||} 
\hline
& Discretisation & $n_{\Gamma,1}$ & $L^2$ error (order)& $H^1$ semi-norm error (order) \\ [0.5ex] 
 \hline\hline
 \multirow{10}{*}{Test 1}&
 \multirow{5}{*}{Conformal}
& 10 & $4.63\times10^{-2} \ $ \ \ \ \ \ \ \ \ & $1.44\times10^{-1} \ $ \ \ \ \ \ \ \ \   \\ 
  \cline{3-5}
& &20 & $1.21\times10^{-2} \ $  (1.93) & $7.09\times10^{-2} \ $ (1.02)   \\ 
 \cline{3-5}
& &40 & $2.82\times10^{-3} \ $ (2.11)   &  $3.45\times10^{-2} \ $ (1.04)  \\ 
 \cline{3-5}
& &80 & $6.99\times10^{-4} \ $  (2.01)& $1.76\times10^{-2} \ $ (0.97)  \\ 
 \cline{3-5}
& &160 & $1.58\times10^{-4} \ $ (2.15)& $9.34\times10^{-3} \ $ (0.90)   \\
 \cline{2-5}\vspace{-1em}
 \\
 \cline{2-5}
& \multirow{5}{*}{Nonconformal}&
10 & $4.61\times10^{-2} \ $\ \ \ \ \ \ \ \  &  $ 1.36\times10^{-1} \ $ \ \ \ \ \ \ \ \  \\ 
 \cline{3-5}
& &20 & $1.25\times10^{-2} \ $  (1.88)  &  $ 7.03\times10^{-2} \ $ (0.95)  \\ 
 \cline{3-5}
& &40 & $2.96\times10^{-3} \ $ (2.08)  & $ 3.53\times10^{-2} \ $ (0.99)  \\ 
 \cline{3-5}
& &80 &  $7.52\times10^{-4} \ $ (1.98) & $ 1.81\times10^{-2} \ $ (0.96)  \\ 
\cline{3-5}
& &160 & $1.66\times10^{-4} \ $  (2.17) & $ 9.46\times10^{-3} \ $ (0.94)  \\
 \hline\hline
 \multirow{10}{*}{Test 2}&
 \multirow{5}{*}{Conformal}
& 10 & $3.79\times10^{-2} \ $ \ \ \ \ \ \ \ \ & $1.69\times10^{-1} \ $ \ \ \ \ \ \ \ \   \\ 
  \cline{3-5}
& &20 & $1.05\times10^{-2} \ $  (1.86) & $8.83\times10^{-2} \ $ (0.94)    \\ 
 \cline{3-5}
& &40 & $2.56\times10^{-3} \ $ (2.03)   &  $4.52\times10^{-2} \ $  (0.96)   \\ 
 \cline{3-5}
& &80 & $6.42\times10^{-4} \ $  (2.00)& $2.34\times10^{-2} \ $ ( 0.95)   \\ 
 \cline{3-5}
& &160 & $1.52\times10^{-4} \ $ (2.08)& $1.28\times10^{-2} \ $ (0.87)   \\ 
 \cline{2-5}\vspace{-1em}
 \\
 \cline{2-5}
& \multirow{5}{*}{Nonconformal}&
10 & $4.32\times10^{-2} \ $\ \ \ \ \ \ \ \  &  $ 1.85\times10^{-1} \ $ \ \ \ \ \ \ \ \   \\ 
 \cline{3-5}
& &20 & $1.09\times10^{-2} \ $  (1.99)  &  $ 8.93\times10^{-2} \ $ (1.05)   \\ 
 \cline{3-5}
& &40 & $2.75\times10^{-3} \ $ (1.98)  & $ 4.67\times10^{-2} \ $ (0.94)   \\ 
 \cline{3-5}
& &80 &  $7.18\times10^{-4} \ $ (1.94) & $ 2.48\times10^{-2} \ $ (0.91)   \\ 
\cline{3-5}
& &160 & $1.67\times10^{-4} \ $  (2.10) & $ 1.32\times10^{-2} \ $ (0.91)   \\ 
 \hline
\end{tabular}
\label{tab:ConvergenceOrder}
\end{table}

We calculate the error in both the $L^2$ norm and the $H^1$ semi-norm (Table \ref{tab:ConvergenceOrder}). The estimated accuracies in the $L^2$ norm and the $H^1$ semi-norm between the reference solution and the exact solution for Test 1 are $2.52\times10^{-5}$ and $1.50\times10^{-3},$ respectively, and for Test 2 are $2.44\times10^{-5}$ and $2.05\times10^{-3},$ respectively. 

In the $H^1$ semi-norm, the method has order 1 convergence, in agreement with the theory of Section \ref{sec:errorAnalysis}. Although we do not have theoretical error estimates for the $L^2$ norm, the method displays numerical convergence of order 2, one order better than in the $H^1$ semi-norm. These results hold irrespective of the discretisation being conformal or nonconformal, with the method over the nonconformal discretisation performing just as well as the method over the conformal one. Albeit, with the nonconformal discretisation, the error is slightly larger than with the conformal one, showing that the consistency error indeed leads to larger overall errors. Additionally, without any special consideration on the dual variable function approximation near the vertices of the polygonal domain, the method is accurate with optimal convergence order in the $H^1$ norm. 
\subsection{Validation via Comparison}\label{sec:2dNumericalValidation}
In \cite{MacDonald:2021:MathBiosci}, we numerically validated a finite difference method for the spatially 1-dimensional system. We now numerically validate our finite element method by comparison with this finite difference scheme. Using the technique of Ludwig et al. \cite{Ludwig:1979:JMathBiol}, we can write the 1-dimensional system on the reference frame as
\begin{subequations}\label{eq:steadyState1DHalfLudwig}
\begin{align}
&\partial_t w_0 = d_0 \partial_{xx} w_0 + c \partial w_0 + w(r - aw) , &  0 < x < L,
\\
&\partial_t w_1 = d_1 \partial_{xx} w_1 + c \partial w_1 - m_1w_1 , & x < 0,
\\ \nopagebreak[4]
&w_0 = \kappa w_1,  & x = 0,
\\
&d_0\partial_{x}w_0 + cw_0 = d_1\partial_x w_1 + c w_1,  &x = 0,
\\ 
&d_0\partial_x w + c w = b w, &x = L,
\\ 
&\lim_{x \to -\infty} w(x, t) = 0.
\end{align}
\end{subequations}
That is, we can rewrite the system on a semi-bounded domain where the boundary coefficient $b = b(\beta, d_2, m_2, c),$ at $x = L$ is a constant \cite{MacDonald:2018:JMathBiol}. Here, $\beta$ is the probability with which an individual moves into the suitable habitat when at $x = L,$ and $d_2,$ and $m_2,$ are the scaled diffusion and mortality rates, respectively for $x > L.$ Steady-state solutions to this system over $x < 0$ are exponential curves; i.e., $w_1 \sim e^{n^+x,}$ where $n^+ > 0,$ and depend on the parameters of the system \cite{MacDonald:2018:JMathBiol}. Thus, over $x < 0,$ far enough away from $x = 0,$ solutions are flat.
\begin{figure}
\centering
\includegraphics[width = 6cm]{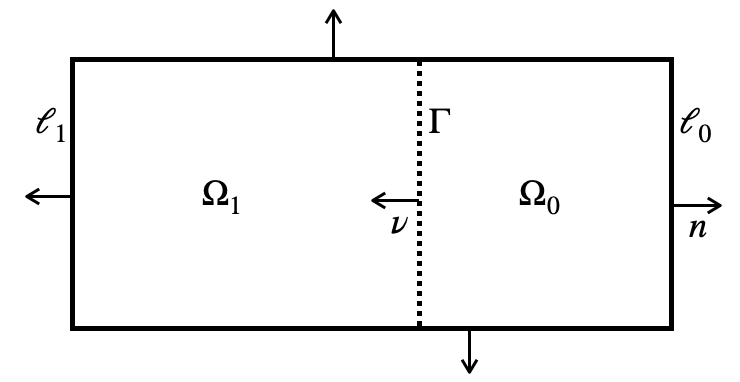}
\caption[One spatial dimension extended to two spatial dimensions]{An illustration of the bi-domain for numerical validation}
\label{fig:biDomIllustrate}
\end{figure}

An analogous 2-dimensional system is
\begin{subequations}\label{eq:biDomain2D}
\begin{align}
&\partial_tw = d_i\Delta w + (\bc\cdot\nabla)w + G(w), &\text{ on }\Omega\backslash\Gamma,
\\
&w_0 = \kappa w_1,  &\text{ on }\Gamma,
\\ 
&d_0\partial_{\bnu}w_0 + (\bc\cdot\nu)w_0 = d_1\partial_{\bnu}w_1 + (\bc\cdot\nu)w_1,  &\text{ on }\Gamma,
\\
&d_0\partial_{\bn} w + (\bc\cdot \bn)w = b w, &\text{ on } \ell_0,
\\
& \partial_nw = 0, &\text{ on } \partial\Omega\backslash\bar{\ell}_0,
\\
&\lim_{\bx \to -\infty} w = 0,
\end{align}
\end{subequations}
where $\Omega = \Omega_1 \cup \Omega_0 \cup \Gamma,$ and for some real values $y_0 < y_1,$ 
\begin{align*}
&\Omega_0 = \{(x, y) \in \mathbb{R}^2 \ | \ 0 < x < L, \ y_0 < y < y_1\}, 
\\
&\Omega_1 = \{(x, y) \in \mathbb{R}^2 \ | \  x < 0, \ y_0 < y < y_1\}, 
\\
&\Gamma = \{(x, y) \in \mathbb{R}^2 \ | \ x = 0, \ y_0 < y < y_1\},
\intertext{and}
&\ell_0 =  \{(x, y) \in \mathbb{R}^2 \ | \ x = L, \ y_0 < y < y_1\}.
\end{align*}
Here, $\bn$ is the outward-pointing unit normal on $\partial\Omega.$ In the special case, $\bc = (c_1, 0),$ this system is spatially homogeneous in the $y$-direction and heterogeneous in the $x$-direction. Cross sections of the steady-state solution of System (\ref{eq:biDomain2D}) along the $x$-direction are the steady-state solution to System (\ref{eq:steadyState1DHalfLudwig}). 

We apply the hybrid finite element method to System (\ref{eq:biDomain2D}) by truncating the subdomain $\Omega_1$ at $\ell_1 = \{(x, y) \in \mathbb{R}^2 |  x = -L', y_0 < y < y_1\},$ where $L' > 0$ (Figure \ref{fig:biDomIllustrate}). On this boundary, we set homogeneous Dirichlet boundary conditions.  The value of $L'$ is chosen large enough so that the solutions are flat near the boundary $\ell_1.$ 

\begin{figure}
\centering
\begin{subfigure}[b]{.6\linewidth}
\includegraphics[scale = 0.4]{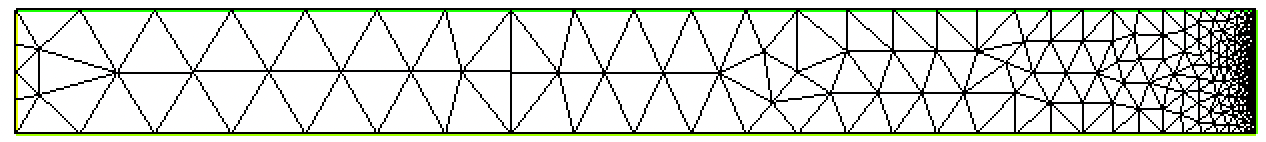}
\caption{$\Omega_1$}
\label{fig:omega1Coarse}
\end{subfigure}
\\
\begin{subfigure}[b]{.35\linewidth}
\includegraphics[scale = 0.25]{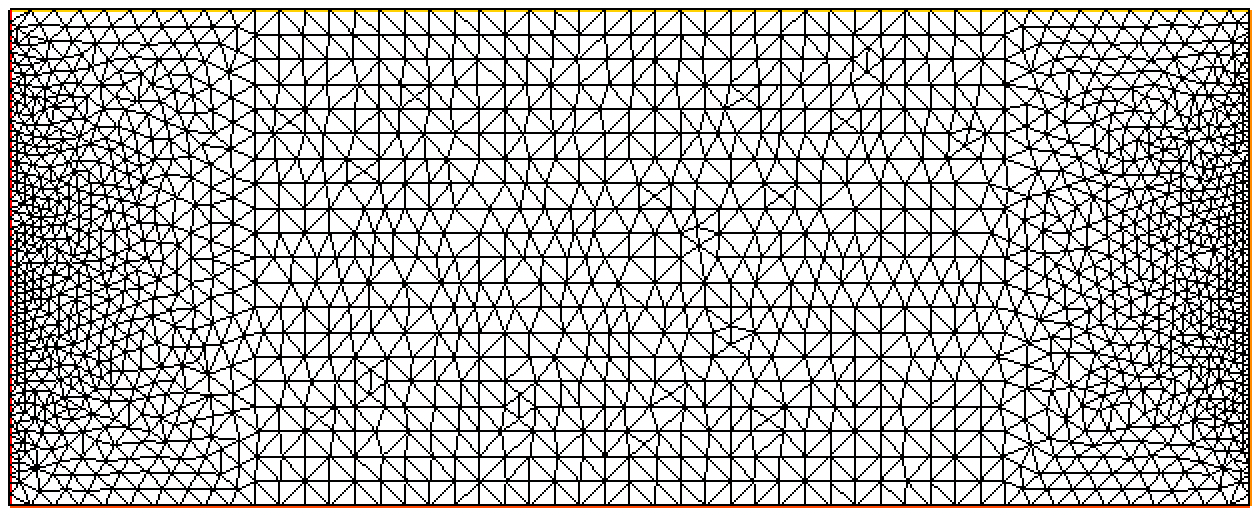}
\caption{$\Omega_0$}
\label{fig:omega0Coarse}
\end{subfigure}
\\
\begin{subfigure}[b]{\linewidth}
\centering
\includegraphics[scale = 0.6]{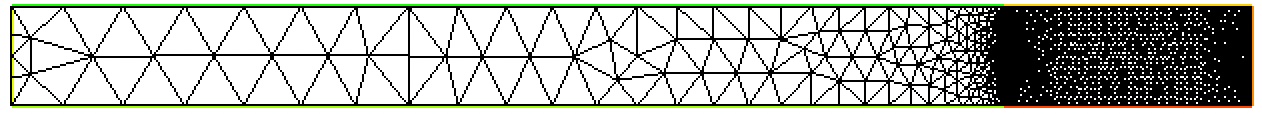}
\caption{$\Omega$}
\label{fig:omegaCoarse}
\end{subfigure}
\caption[The coarsest mesh]{The coarsest triangulation for validation. Panel (a): $\mathcal{T}_h^1$. Panel (b): $\mathcal{T}_h^0,$ stretched for visual ease. Panel (c): $\mathcal{T}_h^1\cup \mathcal{T}_h^0$ }
\label{fig:coarsestGrid}
\end{figure}

We take $\mathbb{P}_1$ approximations for both the Lagrange multiplier and the function $w.$ Inside $\Omega_1$ far away from $\Gamma,$ the solution is nearly constant, so, we use an increasing element size progressing away from $\Gamma$ (Figure \ref{fig:coarsestGrid}). The discretisation is conformal inside each $\Omega_i,$ but across $\Gamma,$ it is not. Thus, we have a nonconformal discretisation. The details of the mesh description are given on a case-by-case basis.

At steady-state, $\partial_tw = 0.$ We let $w^{n}$ denote the numerical solution at time $\tau n,$ where $\tau$ is the uniform time step. We stop the time-stepping iterations when
$$||w^n||_{\tau, 0, \Omega}^2 = \left\Vert\frac{w^{n+1} - w^{n}}{\tau}\right\Vert_{0, \Omega}^2  < 10^{-5}.$$

We denote by $U_{FD}$ the solution of the 1-dimensional model, computed using the finite difference scheme presented in \cite{MacDonald:2021:MathBiosci}. The computational domain is taken as a cut along the $x$-direction of the truncated 2-dimensional domain. We describe the computational grid by the pair $(h_0, \tau_{FD}),$ where $h_0$ is the uniform spacing in the interval $(0, L)$ and $\tau_{FD}$ is the time step. In the interval $(-L', 0),$ we space the grid geometrically with the smallest cell width $h = h_0$ occurring next to the interface $x = 0$ and geometrically increasing away from the interface with a common ratio of $r_c = 1.005.$ For each test case, $U_{FD}$ is determined via mesh refinement and is a grid independent solution within 1\%. For each simulation, we stop the scheme once the forward difference approximation of the time derivative in the $L^\infty$ norm is below $10^{-3}.$ The details of the numerical parameters are given on a case-by-case basis. 

\begin{figure}
\centering
\includegraphics[scale = 0.6]{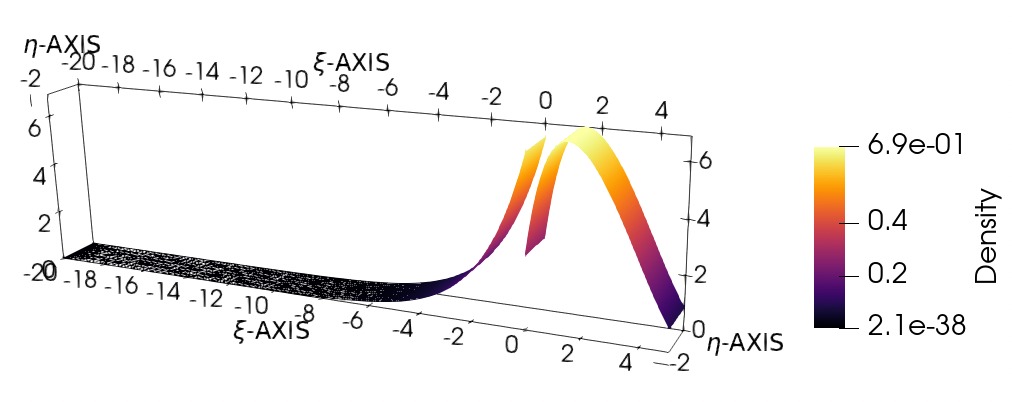}
\caption[A 2-dimensional travelling pulse profile]{A 2-dimensional travelling pulse profile.}
\label{fig:uniformInYDirection}
\end{figure}
The numerical solution obtained by the finite element method is a travelling pulse that is uniform in the $y$-direction (Figure \ref{fig:uniformInYDirection}). Using the open-source software ParaView \cite{Ayachit:2015:Paraview}, we extract the solution along a cut in the $x$-direction, which we denote by $U_{h, \tau},$ where $h$ is a measure of the grid size, and $\tau$ is the time step used in the scheme. We compare the solution $U_{h, \tau}$ to $U_{FD}$ in MATLAB$^{\text{\textregistered}}.$ Here, $U_{h, \tau}$ and $U_{FD}$ are vectors whose $j^{th}$ component is the approximation of the value of the travelling pulse solution over the reference frame at $(\xi_j, \eta_0),$ where $\eta_0$ is a constant. We interpolate $U_{FD}$ onto the extracted grid of $U_{h, \tau}.$ We calculate the relative error, $e_\infty\left(U_{h, \tau}, U_{FD}\right),$ in the infinity norm, where
$$
e_\infty\left(U, V\right) = \frac{\left\| U - V \right\|_{\infty}}{\left\| V \right\|_{\infty}}. 
$$

We choose three different sets of model parameters, called humped shaped, decreasing, and sharply decreasing, which describe the shape of the profile inside the suitable habitat (Table \ref{tab:parameterSetsForValidation}). Since the finite-dimensional scheme is written for the scaled system, we always take $d_0 = r = a = 1.$ 
\begin{table}
\centering
\caption[Parameter sets for validation]{The parameter sets for validation}
\begin{tabular}{| c || c | c | c | c | c | c | c ||} 
 \hline
Test Case & $\alpha$ & $\beta$ & $d_1$ & $d_2$ & $m_1$ & $m_2$ & $c $ \\ [0.5ex] 
 \hline\hline
 \textbf{ humped shaped} & 0.3 & 0.3 & 1 & 1 & 1 & 1 & 1 \\ 
 \hline
 \textbf{decreasing} & 0.8 & 0.5 & 1 &1 & 1 & 0.5 & 1.5   \\ 
 \hline
 \textbf{sharply decreasing} & 0.8 & 0.6 & 2 & 1.3 & 0.1 & 1.4 & 2.5  \\ 
 \hline
\end{tabular}
\label{tab:parameterSetsForValidation}
\end{table}
\begin{figure}
\centering
\begin{subfigure}[b]{.45\linewidth}
\includegraphics[scale = 0.15]{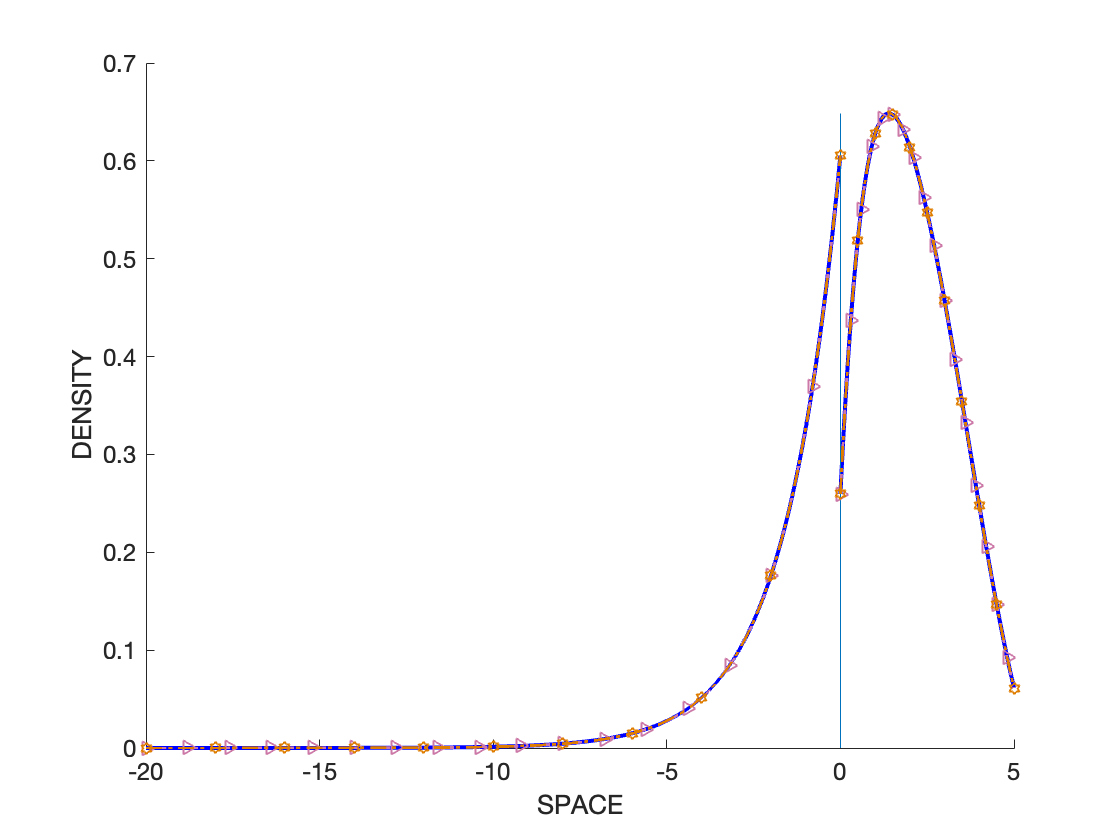}
\caption{Hump shaped}
\label{fig:PS1_P1ApproxFull}
\end{subfigure}
\\
\begin{subfigure}[b]{.45\linewidth}
\includegraphics[scale = 0.15]{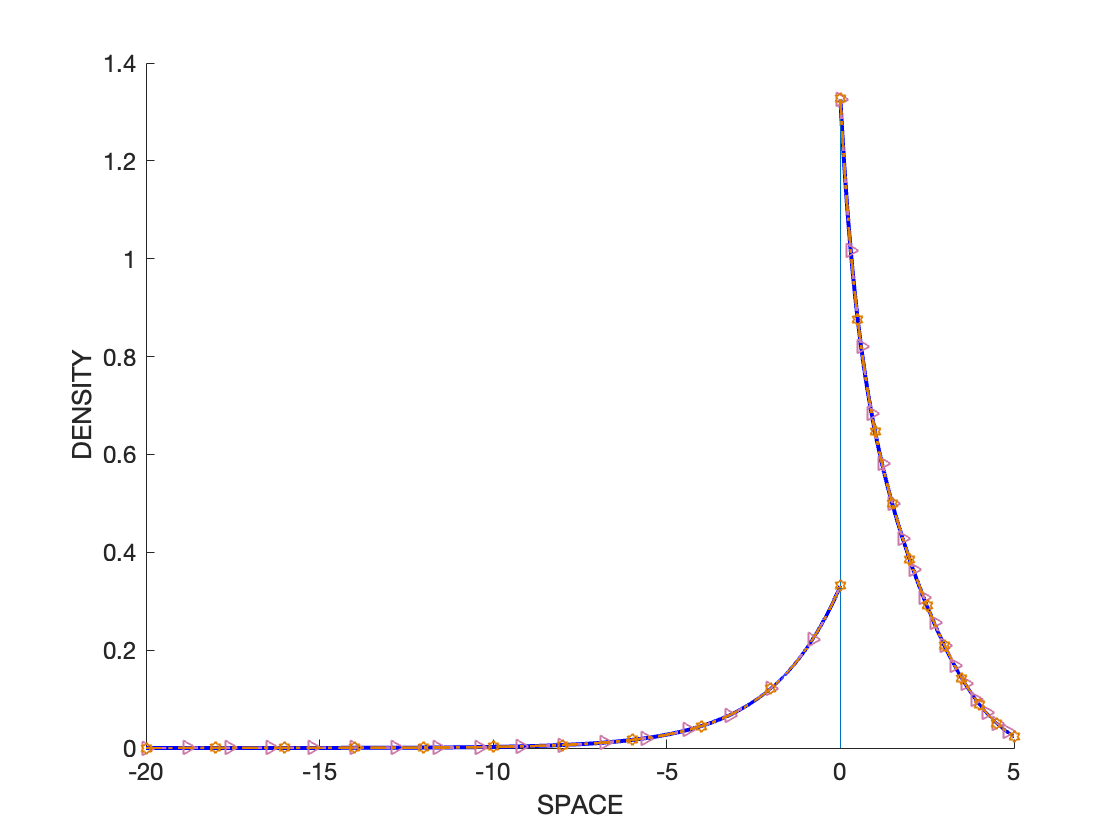}
\caption{Decreasing}
\label{fig:PS2_P1ApproxFull}
\end{subfigure}
\hfil
\begin{subfigure}[b]{.45\linewidth}
\includegraphics[scale = 0.15]{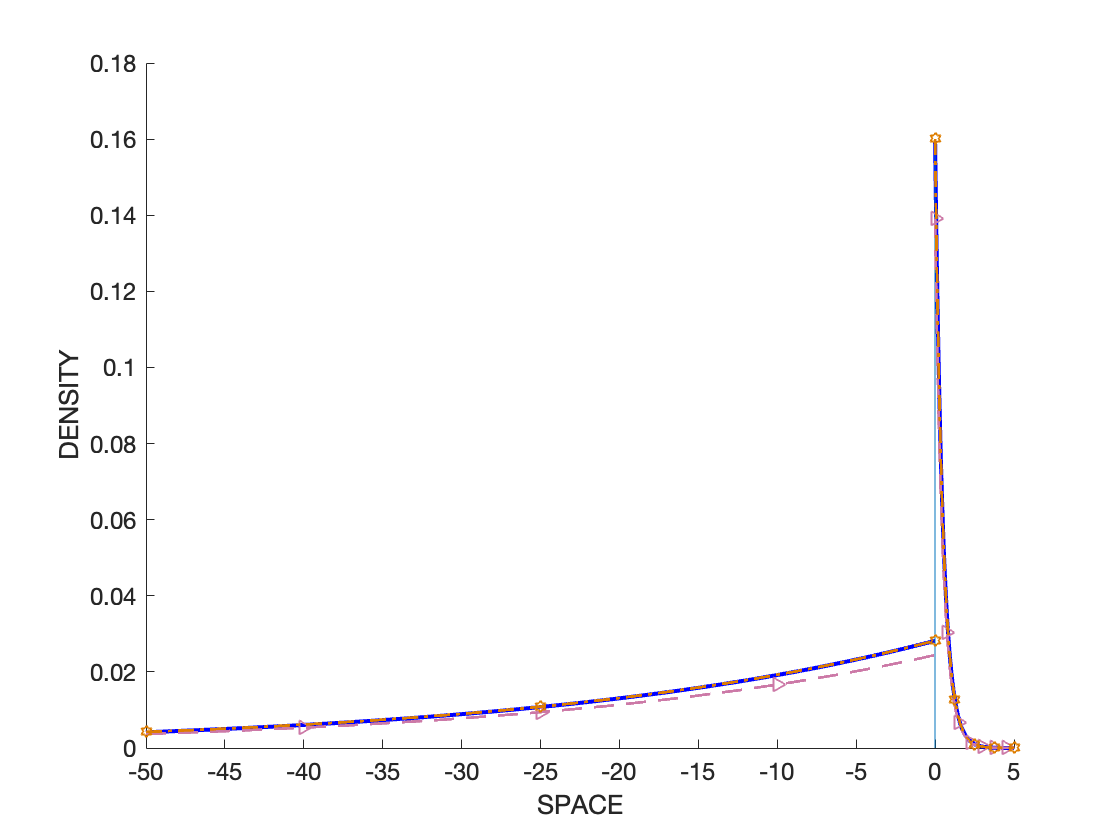}
\caption{Sharply decreasing}
\label{fig:PS3_P1ApproxFull}
\end{subfigure}
\caption[]{Cross sections along the $x$-direction of the 2-dimensional numerical solutions of the finite element method are superimposed on the 1-dimensional numerical solution of the finite difference method. $U_{h, 0.1}$ (dashed with triangle markers, purple), and $U_{\frac{h}{2}, 0.05}$ (dot-dashed with hexagram markers, orange), plotted against $U_{FD}$ (solid, blue). The vertical line is the interface $\Gamma$.}
\label{fig:superImposed}
\end{figure}
\begin{description}
\item[\textbf{humped shaped}]  We run the simulation twice. We take $L = 5,$ and  $L' = 20.$ Our coarsest grid has 30 nodes across the $x$-direction in $\Omega_1,$ and 50 uniformly spaced nodes across the $x$-direction in $\Omega_0$. In the $y$-direction, we place 5 nodes on $\ell_1,$ 50 nodes along both sides of $\Gamma$ and 50 nodes along $\ell_0$ (Figure \ref{fig:coarsestGrid}). Let $h$ be the measure of the largest triangle on our coarsest grid. 

We compute $U_{FD}$ over the grid $(h_0, \tau_{FD}) = (2.5\times10^{-3}, 2.5\times10^{-2}).$ This solution is grid independent to the solution computed on the finer grid $(h_0, \tau_{FD}) = (1.25\times10^{-3}, 1.25\times10^{-3}),$ with a relative error in the infinity norm of approximately $3.4\times10^{-4}.$

The finite element solutions are superimposed on $U_{FD}$ (Figure \ref{fig:superImposed}). The relative error, $e_{\infty}(U_{h, \tau}, U_{FD})$ is below 0.2\% for each solution (Table \ref{tab:relErrorP1Approx}). 
\item[\textbf{decreasing}]  For the finite element approximation we use the same two meshes as we did for parameter set ``humped shaped''. 

We compute $U_{FD}$ over the grid $(h_0, \tau_{FD}) = (1.25\times10^{-3}, 1.25\times10^{-3}).$ This solution is grid independent to the solution computed on the coarser grid $(h_0, \tau_{FD}) = (2.5\times10^{-3}, 2.5\times10^{-2}),$ with a relative error in the infinity norm of approximately $1.3\times10^{-3}.$

The finite element solutions are superimposed on $U_{FD}$ (Figure \ref{fig:superImposed}). All the solutions satisfy $e_{\infty}(U_{h, \tau}, U_{FD}) \leq 0.1\%$ (Table \ref{tab:relErrorP1Approx}). The relative error $e_{\infty}(U_{h, \tau}, U_{FD})$ remains within the bounds of the accuracy prescribed on $U_{FD}.$
\item[\textbf{sharply decreasing}]  Here, we take $L = 5$ and $L' = 250.$ We impose the exact same number of nodes along each boundary as we did in the two previous test cases. In this case, to achieve an exponential curve along a cut in the unsuitable habitat, we stop the time stepping iterations when $||w^n||_{\tau, 0, \Omega}^2  < 10^{-7}.$

We compute $U_{FD}$ over the grid $(h_0, \tau_{FD}) = (7.8125\times10^{-5}, 7.8125\times10^{-4}).$ This solution is approximately grid independent to the solution computed on the coarser grid $(h_0, \tau_{FD}) = (1.5625\times10^{-4}, 1.5625\times10^{-3}),$ with a relative error in the infinity norm of approximately $1.1\times10^{-2}.$ Simulations for these solutions are stopped when the forward difference approximation of the time derivative in the $L^{\infty}$ norm is below $10^{-5}.$ 

All the solutions have the correct qualitative behaviour and the solutions on the finer grids are superimposed on $U_{FD}$ (Figure \ref{fig:superImposed}). On the coarsest grid, the largest error occurs near the interface where the solution visibly differs from $U_{FD}.$ The relative error on the finer mesh is below 1\%.
\end{description}
\begin{table}
\centering
\caption[Relative error for $\mathbb{P}_1$ approximations of the primal variable]{The relative error between the numerical solutions of the finite element solver and $U_{FD},$ for $\mathcal{P}_1$ approximations of the dual and primal variable}
\begin{tabular}{| c ||c c || c|} 
\hline
\multicolumn{4}{|c||}{$e_\infty\left(U_{h, \tau}, U_{FD}\right)$}   \\
 \hline
 \hline
parameter set & $U_{h, 0.1}$ & $U_{\frac{h}{2}, 0.05}$  & $U_{FD}$ accuracy\\ [0.5ex] 
 \hline\hline
 humped shaped & 0.0014 & 0.0006 &  $3.4\times10^{-4}$\\ 
 \hline
 decreasing & 0.0010 & 0.0009 &  $1.3\times10^{-3}$\\ 
 \hline
sharply decreasing & 0.1266& 0.0041  & $1.1\times10^{-2}$ \\ 
 \hline
\end{tabular}
\label{tab:relErrorP1Approx}
\end{table}
\begin{remark}
If we choose $\mathbb{P}_2$ approximations for $w$, then we cannot choose $\mathbb{P}_2$ approximations for $\lambda,$ as the resulting solution is inaccurate along the interface. In this spatial set-up we have that $\Gamma\cap\partial\Omega \neq \varnothing.$ In the literature, this situation needs to be treated specially, where the approximation of the Lagrange multiplier is of polynomial order $\leq s-1$ on $F$ if a vertex of $F$ lies on $\partial\Omega$ \cite{Bernardi:1993:BookDomDecomp}. Going forward, we do not consider situations where $\Gamma\cap\partial\Omega \neq \varnothing.$ We are interested in situations where the suitable habitat is surrounded by unsuitable habitat in all directions. Thus, $\Gamma\cap\partial\Omega = \varnothing.$
\end{remark}

In this section, we compared our finite element solutions to those of the finite difference scheme. We found in all cases that solutions agreed well. We see that with steep curves, care should be taken in the mesh size and the stopping criteria threshold may need to be taken smaller. 
\subsection{Ecologically motivated examples}\label{sec:simulation}
We present an application of the finite element method to two ecologically motivated settings. These examples serve as demonstrations of how our method performs in scenarios beyond the assumptions of our analysis, namely non-polygonal domains and non- constant shifting velocities along $\Gamma.$

\subsubsection{Non-polygonal domains}
Terrestrial habitat is typically two dimensional, whereas almost all existing models consider only one-dimensional simplifications. One way to interpret these one-dimensional models as representations of 2-dimensional landscapes is to assume that the landscape is homogeneous in the perpendicular direction. However, this is a gross simplification that we do not need in our truly 2-dimensional model. We choose a circular disk so that in every direction, there is spatial heterogeneity. We set 
\begin{align*}
&\Omega_0(t) = \{\bx \in \RR^2 \ : \ (x - c_1 t)^2 + (y - c_2 t)^2 < 2\},
\\
&\Omega_1(t) = \{\bx \in \RR^2 \ : \ (x - c_1 t)^2 + (y - c_2 t)^2 > 2 \},
\intertext{and}
&\Gamma(t) = \{\bx \in \RR^2 \ : \ (x - c_1 t)^2 + (y - c_2 t)^2 = 2 \}.
\end{align*} 
Through the change of variable $\xi = x - c_1 t,$ and $\eta = y - c_2 t,$ in the reference frame, we have $\tilde{\Omega}_i = \Omega_i(0)$ and $\tilde{\Gamma} = \Gamma(0).$ For our computational domain, we truncate $\tilde{\Omega}_1$ at a radius of 10. We impose homogeneous Dirichlet boundary conditions along $\partial\tilde{\Omega}.$ We use a conformal discretisation where we place 80 nodes along $\partial\tilde{\Omega}$ and 160 nodes along $\tilde{\Gamma}.$ We take a mesh size enlarging from $\tilde{\Gamma}$ out towards $\partial\tilde{\Omega}.$ We take a uniform time step of $\tau = 0.025.$ For an initial condition, we take the Gaussian function
$$
w(\xi, \eta, 0) = \frac{40}{\pi}\exp\left(-0.5\left(\left(\frac{\xi}{0.5}\right)^2 + \left(\frac{\eta}{0.5}\right)^2\right)\right). 
$$
Solutions of System (\ref{eq:master2D}) with a constant shift approach a unique travelling pulse solution; see \cite{Berestycki:2008:DiscretContinDynS} and Section \ref{sec:numericalConvergence}. Such a solution is a steady-state solution in the reference frame, and so $\partial_tw = 0.$ As in the previous section, we employ this as a stopping condition and stop the iterations when $||w^n||_{\tau, 0, \Omega}^2 < 10^{-5}.$

Our finite element method was analysed under the assumption that the interface was polygonal. Over this computational domain, $\Gamma$ is a smooth curve. We ran simulations for multiple parameter sets and found that solutions were accurate provided the interface is well represented by the linear mesh. 

\begin{figure}
\centering
\begin{subfigure}[b]{0.45\linewidth}
\centering
\includegraphics[scale = 0.15]{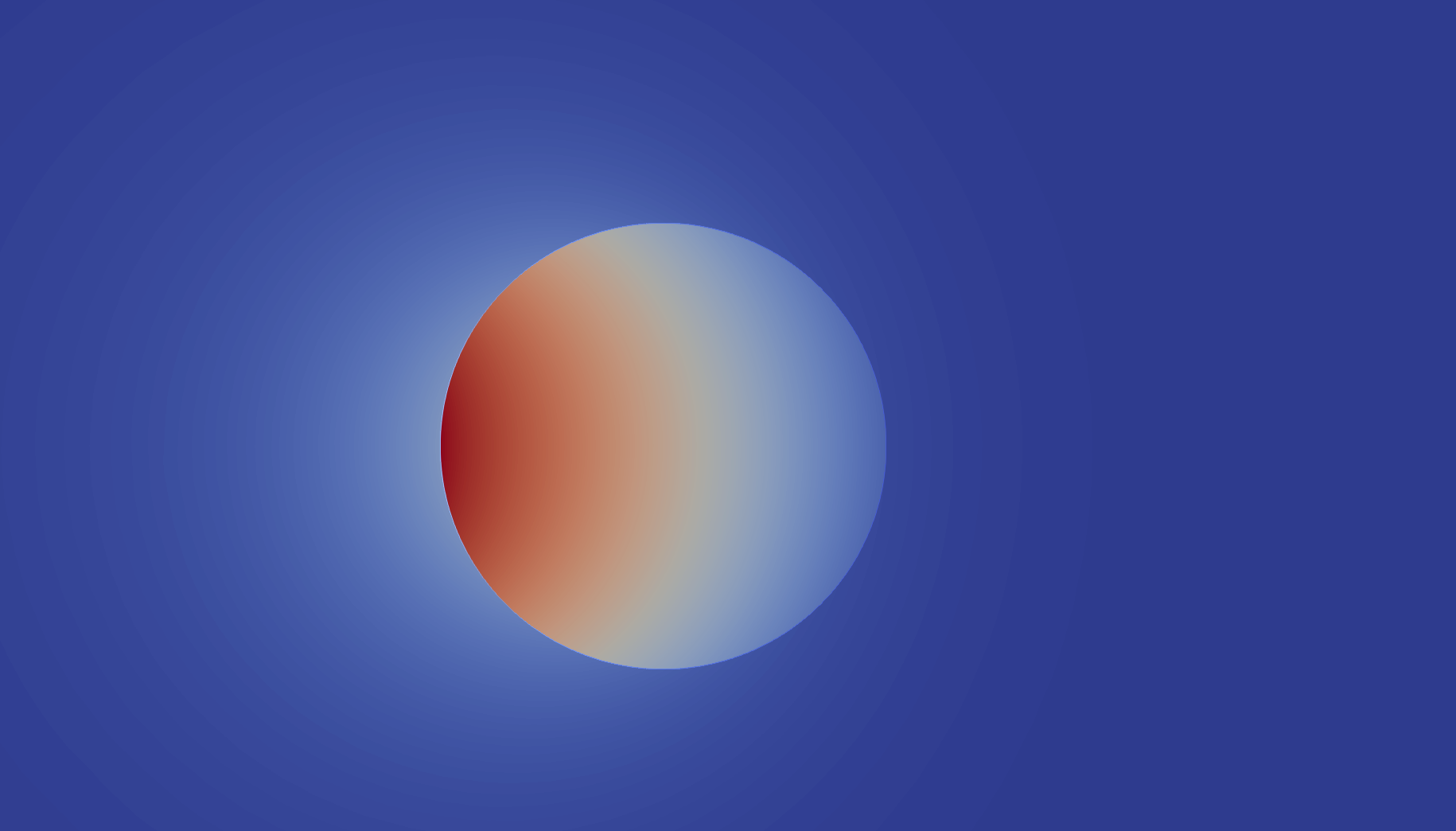}%{2dalpha7D120D010c110c20r12a8m10PLANE.png}
\caption{High bias}
\label{fig:2dalpha7D120D010c110c20r12a8m10PLANE}
\end{subfigure}
\hfil
\begin{subfigure}[b]{0.45\linewidth}
\centering
\includegraphics[scale = 0.15]{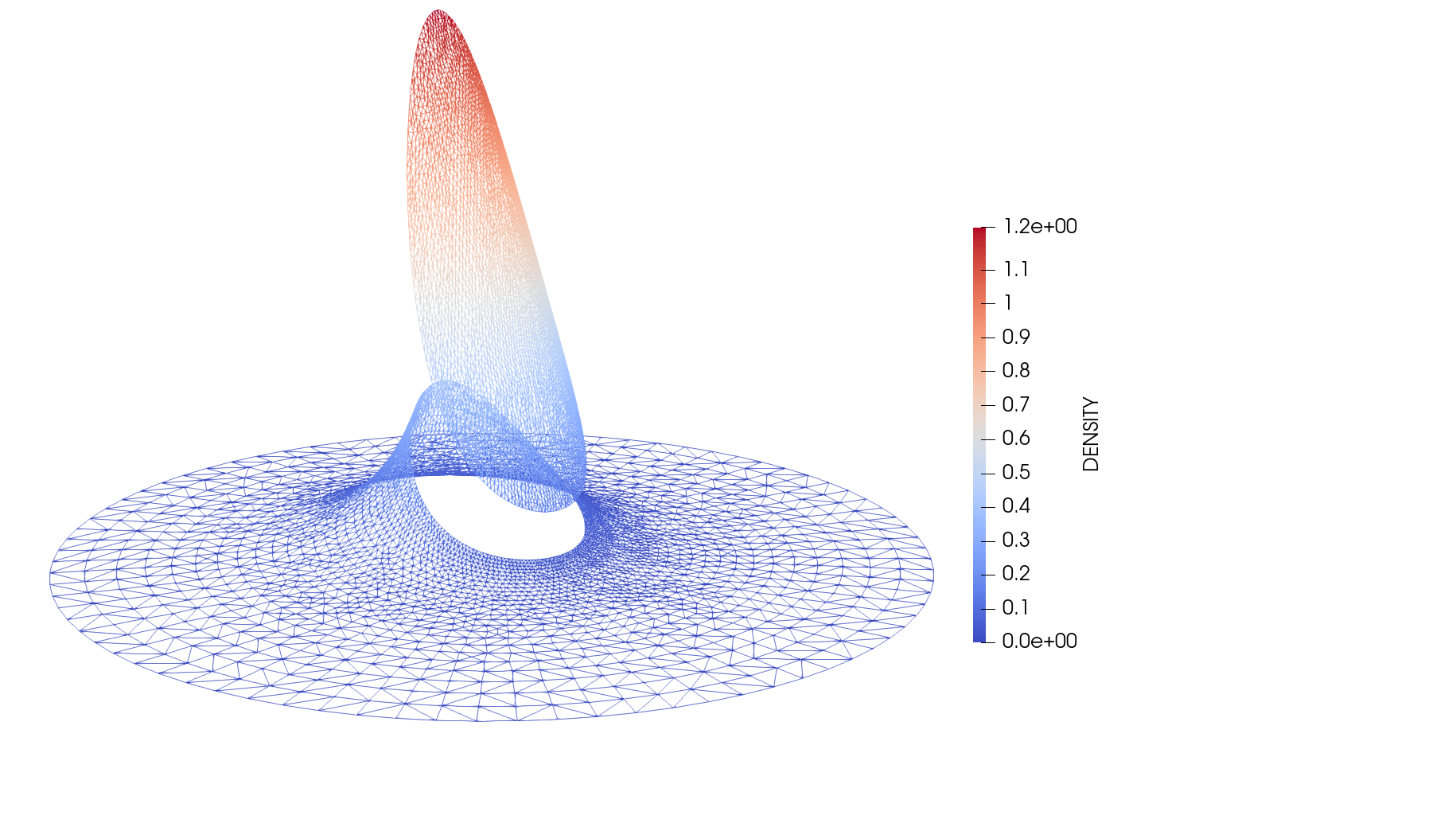}%{2dalpha7D120D010c110c20r12a8m10SURFACE.png}
\caption{High bias}
\label{fig:2dalpha7D120D010c110c20r12a8m10SURFACE}
\end{subfigure}
\begin{subfigure}[b]{0.45\linewidth}
\centering
\includegraphics[scale = 0.15]{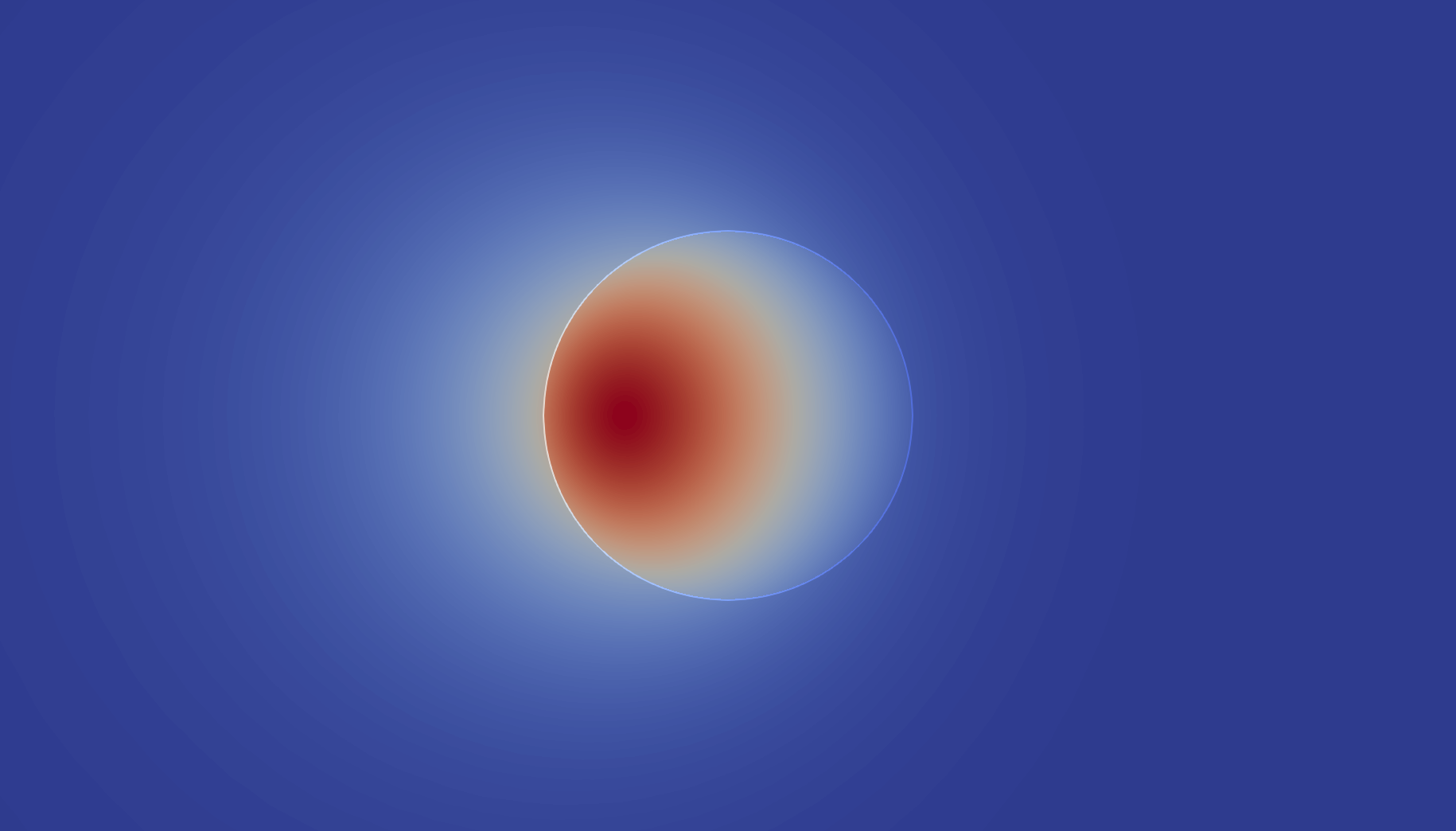}%{2dalpha5D120D010c110c20r12a8m10PLANE.png}
\caption{No bias}
\label{fig:2dalpha5D120D010c110c20r12a8m10PLANE}
\end{subfigure}
\hfil
\begin{subfigure}[b]{0.45\linewidth}
\centering
\includegraphics[scale = 0.15]{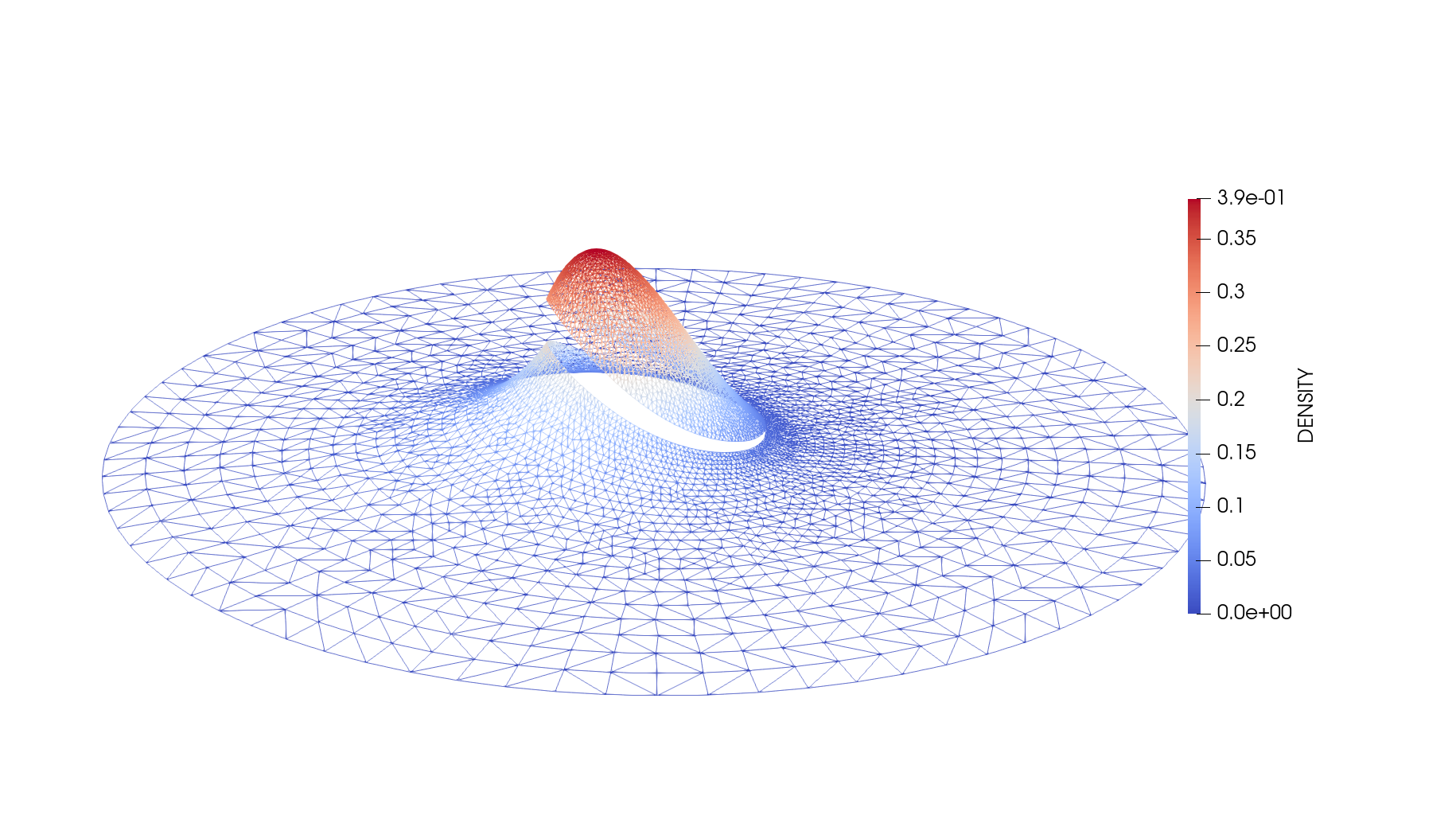}%{2dalpha5D120D010c110c20r12a8m10SURFACE.png}
\caption{No bias}
\label{fig:2dalpha5D120D010c110c20r12a8m10SURFACE}
\end{subfigure}
\caption[Density over a circular suitable habitat]{Density over a circular suitable habitat in the reference frame. Panels (a) and (b): $\alpha = 0.7.$ Panels (c) and (d): $\alpha = 0.5.$ All other parameter values are $d_1 = 2,$ $d_0 = 1,$ $\bc = (1, 0),$ $r = 1.2,$ $a = 8,$ $m = 1.$ Density scales are different across rows of panels.}
\label{fig:circularDomain}
\end{figure}

We present the travelling pulse solution for two of these parameter sets, described in the caption of Figure \ref{fig:circularDomain}. For these cases, we consider the $x$ and $y$ axes as the latitudinal and longitudinal coordinates, respectively. Thus, the poleward environmental shift \cite{Chen:2011:Science} is represented by a constant shift along the positive $x$-direction. The population density profile has a local maximum on the $x$-axis (Figure \ref{fig:circularDomain}). The shift causes the maximum of the density profile within the suitable habitat to skew towards the direction opposite of the shift. The location of the local maximum of the density profile over the suitable habitat depends on model parameters. For example, when there is no bias for the suitable habitat, this local maximum occurs within the suitable habitat (Figures \ref{fig:2dalpha5D120D010c110c20r12a8m10PLANE} and \ref{fig:2dalpha5D120D010c110c20r12a8m10SURFACE}). Increasing the bias, causes the local maximum to shift towards the boundary $\Gamma,$ where a cut of the density profile along the $x$-axis is a decreasing curve (Figures \ref{fig:2dalpha7D120D010c110c20r12a8m10PLANE} and \ref{fig:2dalpha7D120D010c110c20r12a8m10SURFACE}).  
\subsubsection{Non-constant shifting speeds}
\begin{figure}
\centering
\begin{subfigure}[b]{0.45\linewidth}
\centering
\includegraphics[scale = 0.15]{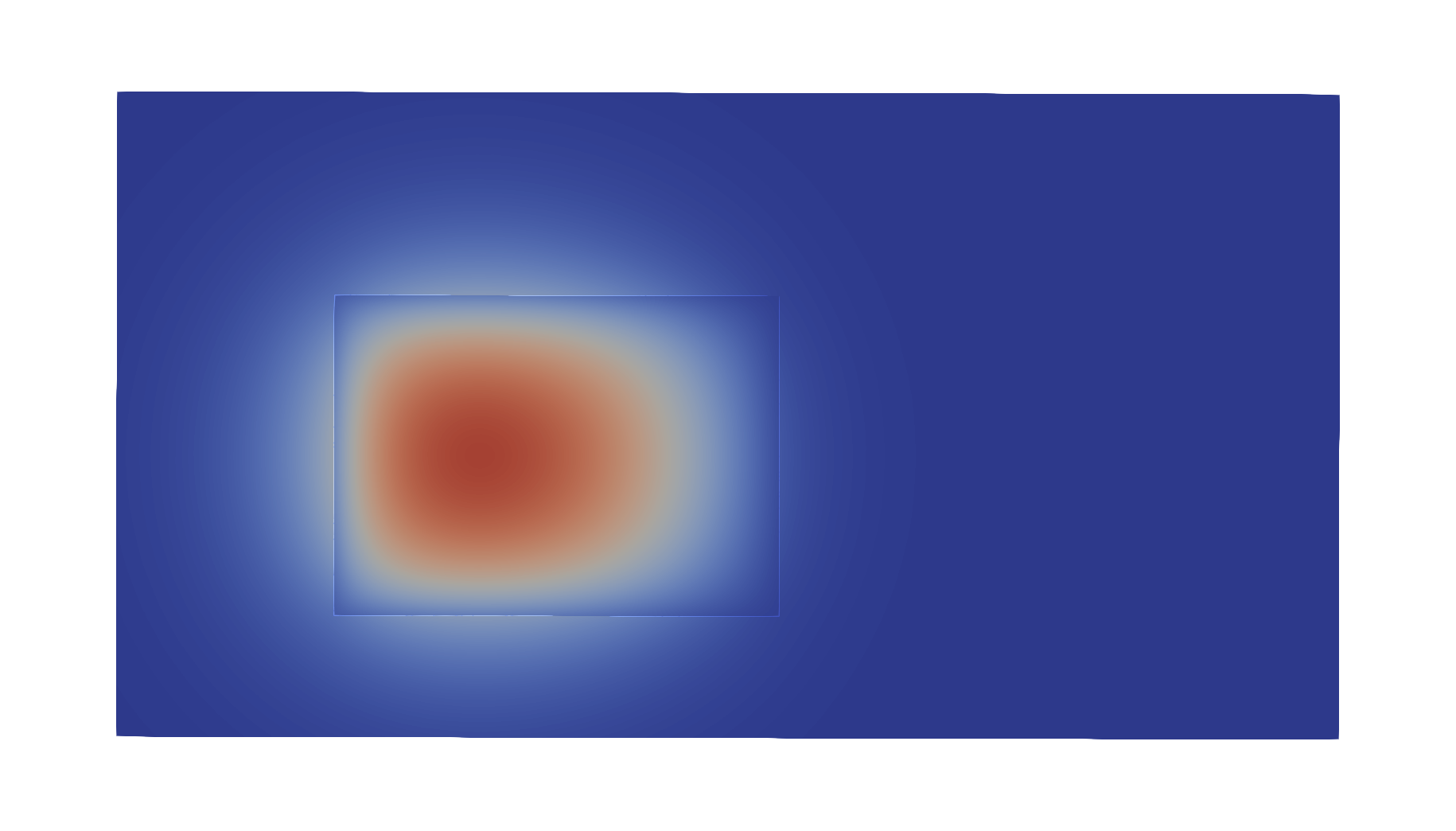}
\caption{$t = 4$}
\label{fig:t4}
\end{subfigure}
\hfil
\begin{subfigure}[b]{0.45\linewidth}
\centering
\includegraphics[scale = 0.15]{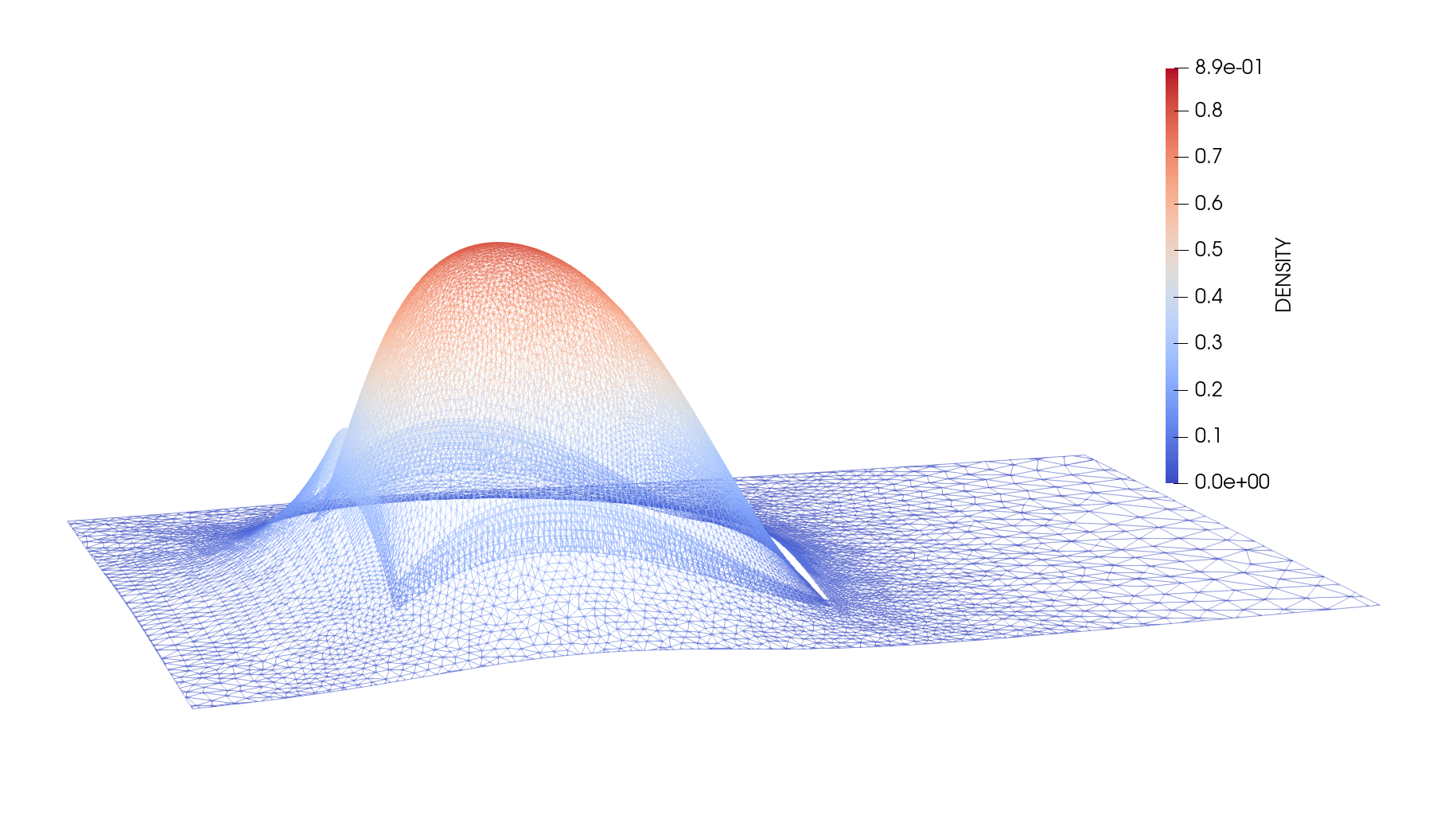}
\caption{$t = 4$}
\label{fig:t4_3D}
\end{subfigure}
\begin{subfigure}[b]{0.45\linewidth}
\centering
\includegraphics[scale = 0.15]{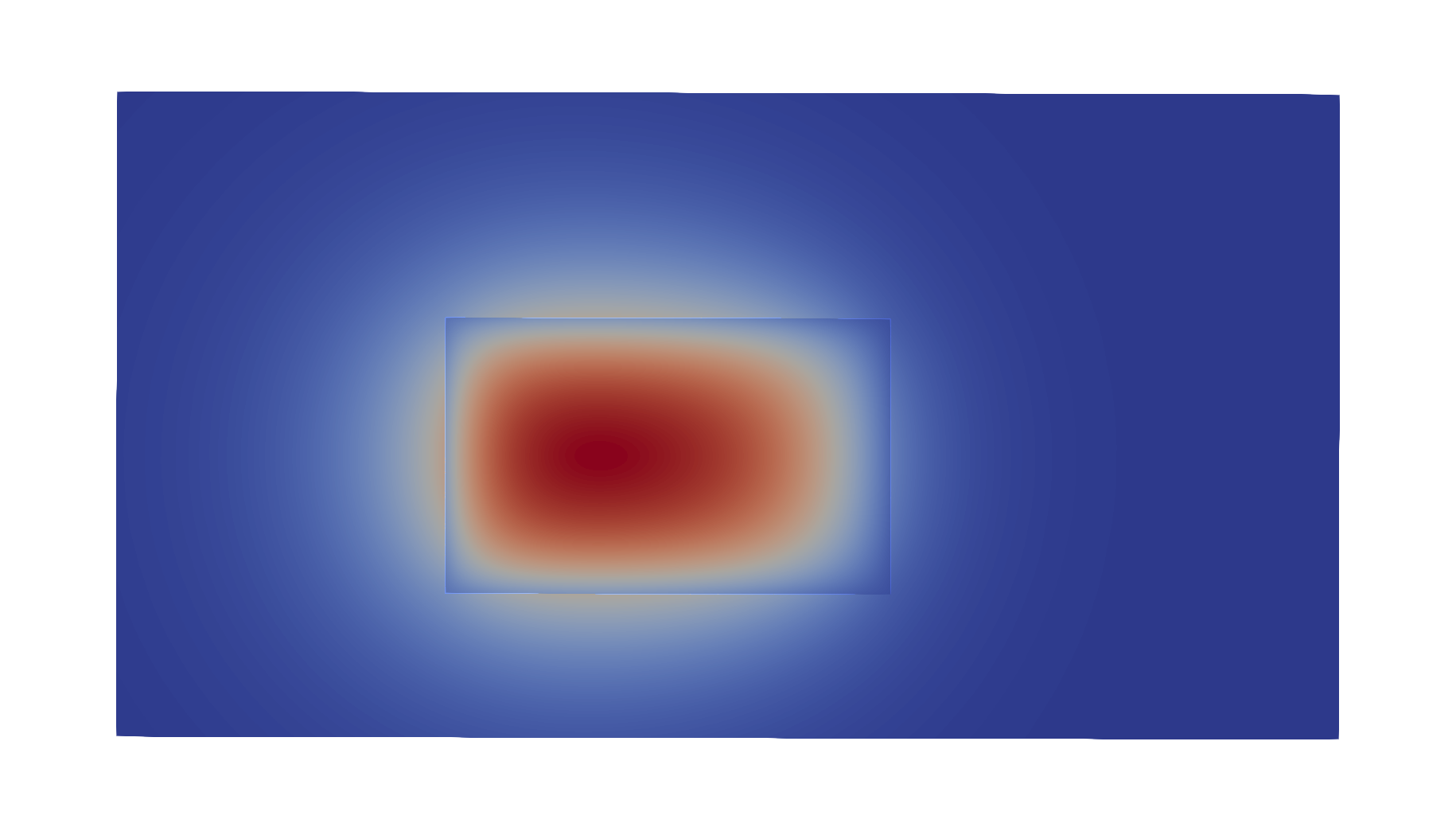}
\caption{$t = 9$}
\label{fig:t9}
\end{subfigure}
\hfil
\begin{subfigure}[b]{0.45\linewidth}
\centering
\includegraphics[scale = 0.15]{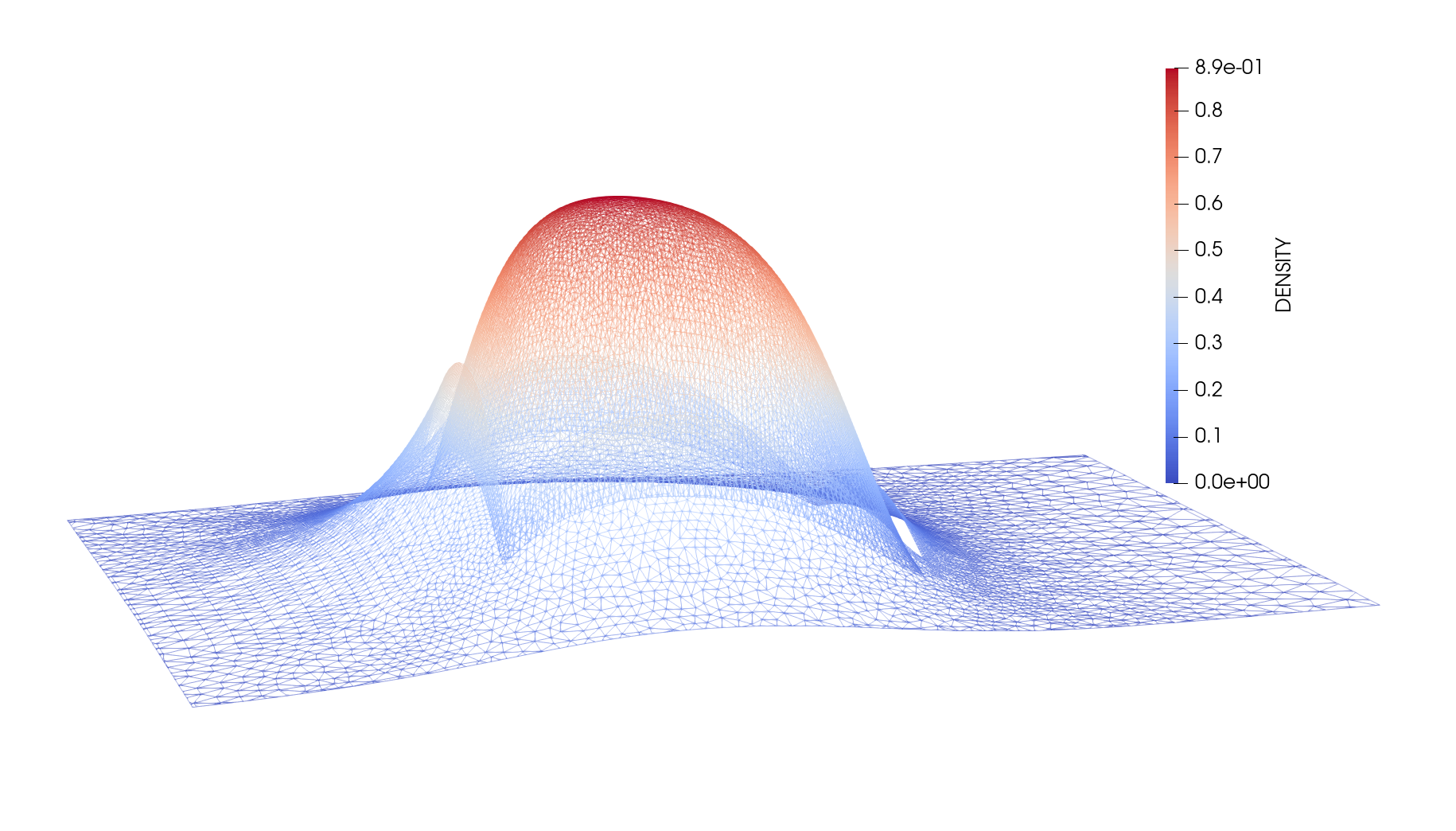}
\caption{$t = 9$}
\label{fig:t9_3D}
\end{subfigure}
\begin{subfigure}[b]{0.45\linewidth}
\centering
\includegraphics[scale = 0.15]{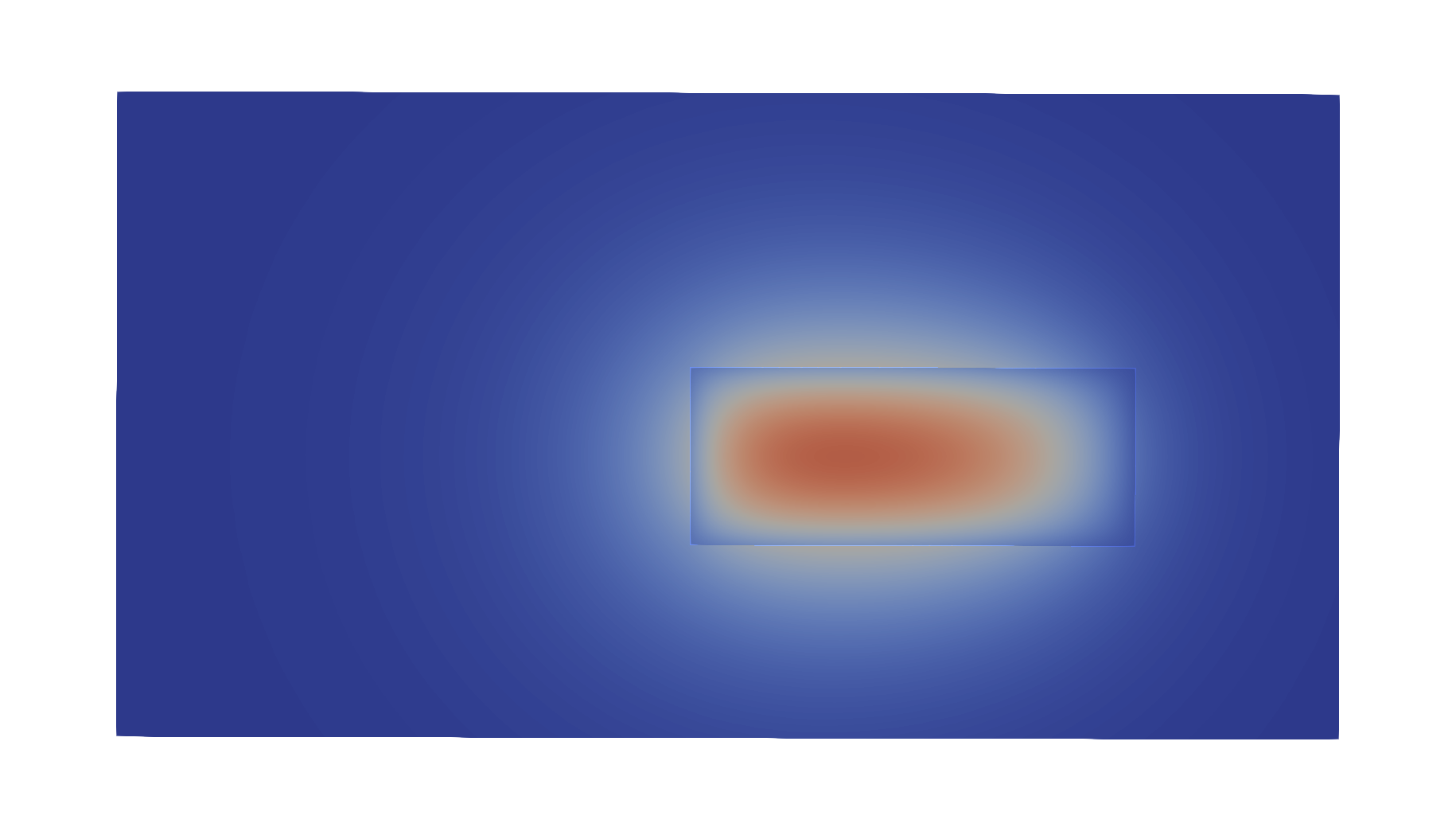}
\caption{$t = 20$}
\label{fig:t20}
\end{subfigure}
\hfil
\begin{subfigure}[b]{0.45\linewidth}
\centering
\includegraphics[scale = 0.15]{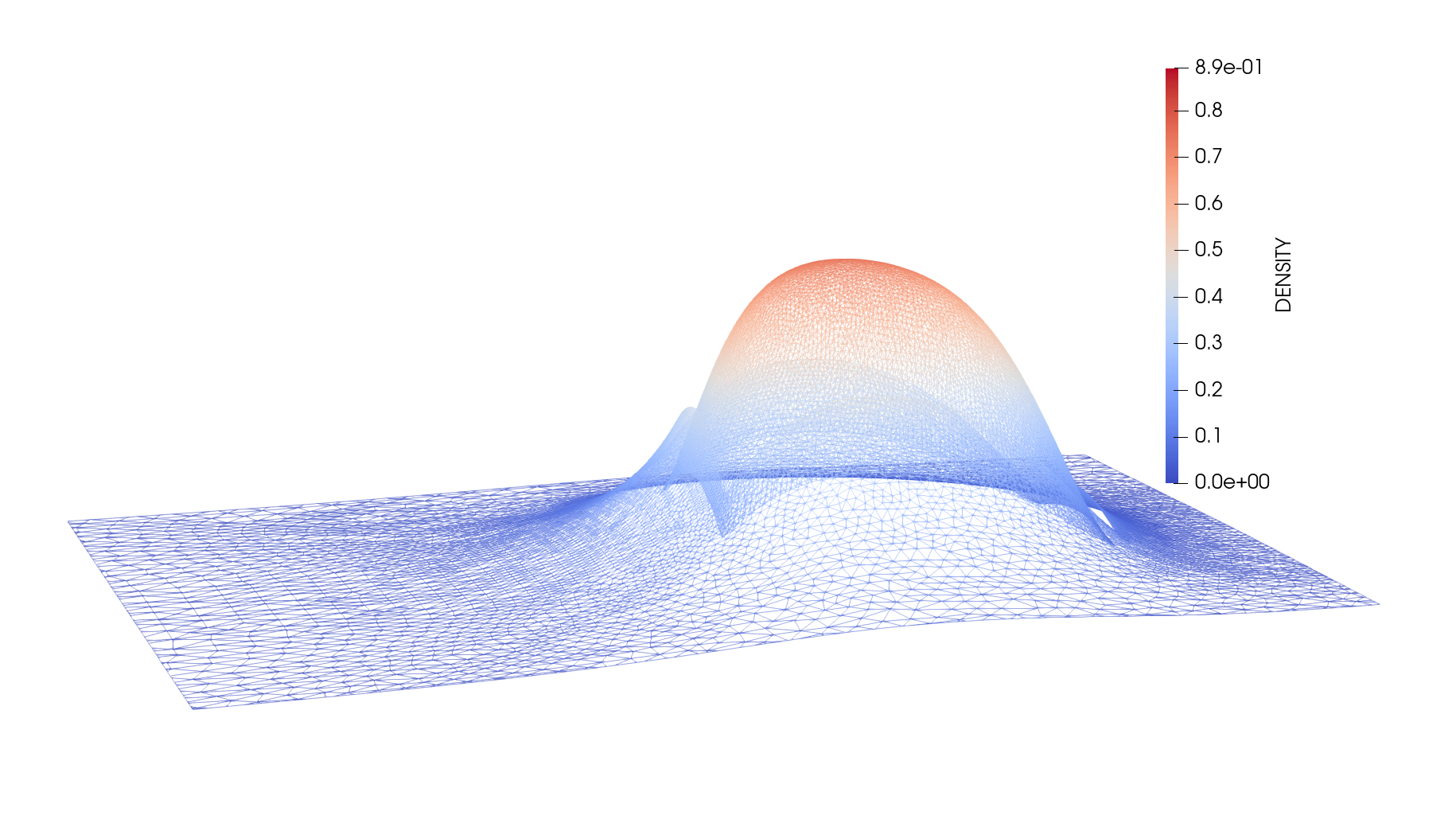}
\caption{$t = 20$}
\label{fig:t20_3D}
\end{subfigure}
\caption[Density over a shrinking habitat]{Density over a shrinking habitat in the physical frame. Parameter values are $d_1 = 2,$ $d_0 = 1,$ $\alpha = 0.3,$ $c_1= 0.5,$ $c_2 = 0.1,$ $r = 1,$ $a = 1,$ and $m = 0.5.$}
\label{fig:shrinkingDomain}
\end{figure}
Change in topography and land connectivity can influence the geometry of the suitable habitat. For example, the suitable habitat may narrow in size as it squeezes between two mountains or shifts onto a peninsula. To capture such a scenario, we take

\begin{align*}
&\Omega_0(t) = \{\bx \in \RR^2 \ : \  c_1t  < x < L_x + c_1t, \ -L_y + c_2t < y < L_y - c_2t\ \},
\\
&\Omega_1(t) = \RR^2 \backslash \bar{\Omega}_0(t),
\intertext{and}
&\Gamma(t) = \cup_{i = 1}^4 \Gamma_i,
\end{align*} 
where 
\begin{align*}
&\Gamma_1(t) = \{\bx \in \RR^2 \ : \  c_1t <  x < L_x + c_1t, \ y = -L_y + c_2t \ \},
\\
&\Gamma_2(t) = \{\bx \in \RR^2 \ : \    x = L _x+ c_1t, \ -L_y + c_2t < y < L_y - c_2t \},
\\
&\Gamma_3(t) = \{\bx \in \RR^2 \ : \  c_1t <  x < L_x + c_1t, \  y = L_y - c_2t   \},
\\
&\Gamma_4(t) = \{\bx \in \RR^2 \ : \    x = c_1t , \ -L_y + c_2t < y < L_y - c_2t  \},
\end{align*} 
where $L_x, \ L_y,$ and $c_i$ are all positive constants. In this case, only transient times up to some $t_f < T = \frac{L_y}{c_2}$ are relevant since the suitable habitat is reduced to a line segment at time $T.$

The transformation to the fixed reference frame is 
$$\xi = x - c_1t, \quad \eta = y\frac{L_y}{L_y - c_2t},$$
so that $\Omega_i(t)$ is mapped to $\Omega_i(0).$

We implement this scenario with $\Omega_0 = [-4, 4] \times [0, 10].$ We use a conformal discretisation and place 80 nodes on each side of the interior rectangle and 240 nodes on each side of the exterior rectangle. Again, we take a mesh size enlarging from $\tilde{\Gamma}$ out towards $\partial\tilde{\Omega}.$ Our time step is $\tau = 0.01.$ For an initial condition we take the Gaussian distribution 
$$w(\xi, \eta, 0) = \frac{20}{\pi}\exp\left(-0.5\left(\left(\frac{\xi-5}{0.5}\right)^2 + \left(\frac{\eta}{0.5}\right)^2\right)\right).$$ 
Through out our analysis in previous sections, we assumed a constant shifting speed $\bc$ along $\Gamma.$ Here $\bc$ has a piecewise construction along $\Gamma.$ Using Equations (\ref{eq:wTransformationGeneral}) and (\ref{eq:fluxRef}) we can formulate the corresponding hybrid formulation for which we apply our finite element method.

Figure \ref{fig:shrinkingDomain} illustrates how the population density changes over a shrinking domain with model parameters given in the caption. With the given initial conditions, the population first grows in that higher local densities are viewed across the suitable habitat within the first nine time units. As time marches forward, the population eventually starts to feel the effect of the shrinking domain and starts to die out as is clearly evident after the first 20 time units. 

\section{Summary}
In the present paper, we developed a 2-dimensional moving-habitat model in a reaction-diffusion framework. The system can capture arbitrary movement of the interface and jumps in density across the interface. The main scope of this work is in developing and analysing a hybrid variational formulation and a corresponding finite element method that can capture a nonstandard jump in density over a shifting interface. For analytical tractability, we assumed that the shifting speed of the interface was constant or, as in the case of the error analysis, zero. Yet, in implementation, we showed that with nonzero constant shifting speeds our numerical method remains optimally convergent. Moreover, in our error analysis, we do not need to impose special conditions on the space of the Lagrange multiplier of our hybrid finite element method to reach an optimal error estimate on the primal
variable, in our case the population density. This is in sharp contrast with the literature on mortar finite element methods where special modifications of the finite element discretisation of the Lagrange are always done near domain corners. These tedious modifications are not needed for the convergence and accuracy of the primal variable. 

The assumptions of our analysis are in general not necessary for implementation of the method. We demonstrated the strength of our method over a non-polygonal domain and over a shrinking domain, both of which are ecologically motivated and inspired by the studies of moving-habitat models. In the spirit of conservation, our method is an optimal tool to identify mechanisms that can help to determine a species' susceptibility to extinction in a wide variety of scenarios for climate-driven moving habitats. With this algorithm in hand, immediate future projects will move towards heavily exploiting the second dimension. For example to gain insight on the impact of the habitat geometry on population persistence (in both moving and stationary domains) and studying population dynamics in transient times, applicable to cases such as an ever shrinking domain, where the asymptotic outcome of the population is inevitably extinction.

\section*{Acknowledgements}
This work was supported by the NSERC CREATE ‘‘Simulation Based Engineering’’ program and an Ontario Graduate Scholarship to JSM, an NSERC, Canada Discovery Grant (RGPIN-2016-0495) to FL and an NSERC, Canada Discovery Grant (RGPIN-2019-06855) to YB.

 \section*{Conflict of Interest}

 The authors declare that they have no conflict of interest.

\appendix

\section{Deriving the flux condition}\label{sec:AppA}

\begin{figure}
\centering
\includegraphics[width=6.5cm,height=6cm,angle=0]{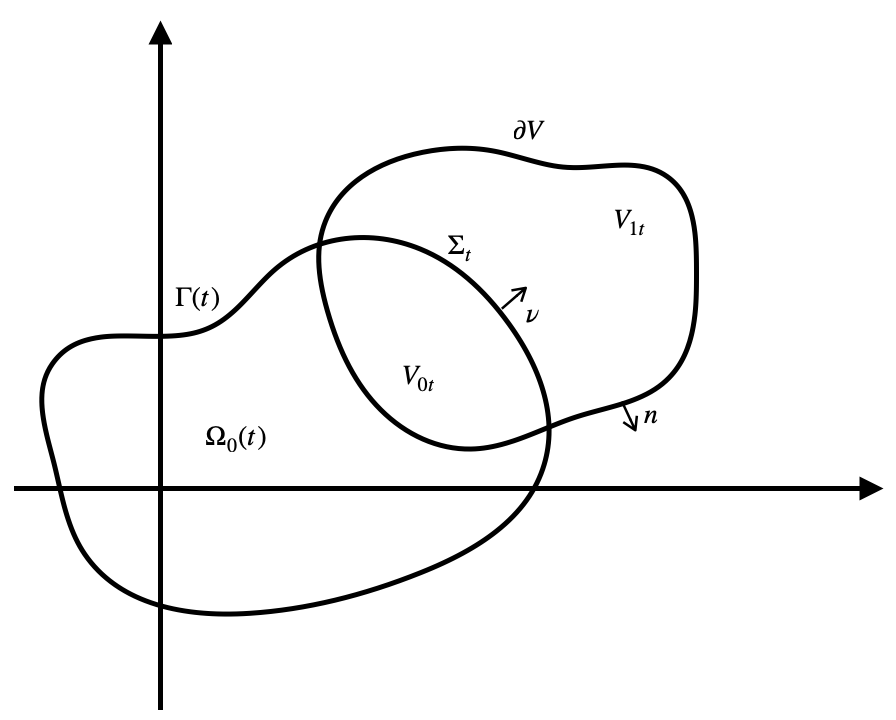}
\caption[Control volume intersecting with the suitable habitat]{A control volume $V$ intersecting with $\Omega_0(t).$}
\label{fig:fluxDerivPicture}
\end{figure}
Here we derive the flux condition Equation (\ref{eq:FluxCondition2D}). Let $V$ denote an arbitrary open control volume (Figure \ref{fig:fluxDerivPicture}). Then, from conservation laws, the change of mass in the control volume $V$ is equal to the diffusive flux across the boundary plus the change in mass due to the source term; i.e., 
\begin{equation}\label{eq:conservationLaw}
\frac{d}{dt}\int_Vu\mathrm{d}\bx = -\int_{\partial V} F(u)\cdot n \mathrm{d}\bx + \int_VG(u)\mathrm{d}\bx,
\end{equation}
where $\partial V$ denotes the boundary of $V,$ $n$ denotes the outward pointing unit normal on $\partial V,$ $F(u) = -D(x,y)\nabla u$ (Fick's Law for diffusive flux), $D(x,y) = d_i,$ on $\Omega_i(t),$ and 
$$G(u) = \begin{cases} u(r - au), &\text{ on } \Omega_0(t), 
\\ -mu, &\text{ on } \Omega_1(t).\end{cases}$$

Let $\Sigma_t = \Gamma(t)\cap V,$ $V_{it} = \Omega_i(t)\cap V,$ for $i = 0, 1.$ From the previous definitions of $\bnu$ and $\Sigma,$ $\bnu$ is also the unit normal pointing from $V_{0t}$ into $V_{1t}.$ Similarly, $c$ is the velocity of $\Sigma_t.$ 

The left-hand side of equation (\ref{eq:conservationLaw}) is then 
\begin{align}
\frac{d}{dt}\int_Vu \ \rd \bx &=  \sum_{i = 0}^1\frac{d}{dt}\int_{V_i}u \ \rd \bx,
\notag\\
&= \sum_{i = 0}^1\int_{V_{it}}u_{it} \ \rd \bx + \int_{\Sigma_t}u_0\bc\cdot \bnu \ \rd s - \int_{\Sigma_t}u_1\bc\cdot \bnu \ \rd s, 
\notag\\
&= \sum_{i = 0}^1\left[\int_{V_{it}}d_i\Delta u_i \ \rd \bx + \int_{V_{it}}G(u_i) \ \rd \bx \right]  
\notag\\ 
&\qquad+ \int_{\Sigma_t}u_0\bc\cdot \bnu \ \rd s- \int_{\Sigma_t}u_1\bc\cdot \bnu \ \rd s,
\notag\\
&= \sum_{i = 0}^1\left[\int_{\partial V_{it} \cap \partial V}d_i\frac{\partial u_i}{\partial n} \ \rd s + \int_{V_{it}}G(u_i) \ \rd \bx\right] +\int_{\Sigma_t}d_0 \partial_{\bnu}u_0 \ \rd s 
\notag\\
&\qquad - \int_{\Sigma_t}d_1\partial_{\bnu}u_1 \ \rd s +\int_{\Sigma_t}u_0\bc\cdot \bnu  \ \rd s - \int_{\Sigma_t}u_1\bc\cdot \bnu  \ \rd s, 
\notag \\
&= -\int_{\partial V} F(u)\cdot n  \ \rd \bx + \int_VG(u) \ \rd \bx +\int_{\Sigma_t}d_0 \partial_{\bnu}u_0 \ \rd s 
\notag \\ \label{eq:controlVolumeCalc}
&\qquad - \int_{\Sigma_t}d_1\partial_{\bnu}u_1 \ \rd s +\int_{\Sigma_t}u_0\bc\cdot \bnu  \ \rd s - \int_{\Sigma_t}u_1\bc\cdot \bnu  \ \rd s. 
\end{align}
The second line is an application of the differentiation formula for moving regions; see Appendix C of \cite{evans:2010:bookPDE}. (In fluid dynamics, this differentiation formula is often called Reynold's Transport Theorem \cite{lidstrom:2011:MathMechSolids}.)
Thus, from equations (\ref{eq:conservationLaw}) and (\ref{eq:controlVolumeCalc})
$$\int_{\Sigma_t}d_0 \partial_{\bnu}u_0 \ \rd s - \int_{\Sigma_t}d_1\partial_{\bnu}u_1 \ \rd s +\int_{\Sigma_t}u_0\bc\cdot \bnu  \ \rd s - \int_{\Sigma_t}u_1\bc\cdot \bnu  \ \rd s = 0.$$
Since $V$ is arbitrary, it follows that
$$d_0\partial_{\bnu}u_0 + (\bc\cdot \bnu)u_0 = d_1\partial_{\bnu}u_1+ (\bc\cdot \bnu)u_1, \quad \bx \in \Gamma(t).$$

One can derive the same equation through similar calculations by stating that mass is conserved over $\mathbb{R}^2$ in the absence of population dynamics, as we previously did in \cite{MacDonald:2018:JMathBiol}. 
%%-----------------------------
%%      your bibliography
\bibliographystyle{abbrv}
\bibliography{hybridMethodMovingHabitatBib.bib}
%%-----------------------------
\end{document}